\documentclass[preprint,12pt]{elsarticle}

\usepackage[utf8]{inputenc} 
\usepackage[T1]{fontenc}    
\usepackage{amsmath,amsthm, amssymb,amsfonts}
\usepackage{mathtools, bm, xcolor}
\usepackage[ruled,linesnumbered,vlined]{algorithm2e}
\usepackage{algorithmic}
\usepackage{hyperref}       
\usepackage{url}            
\usepackage{booktabs}       
\usepackage{nicefrac}       
\usepackage{microtype}      
\usepackage{subcaption}

\newtheorem{defin}  {\hspace*{1pc} Definition}

\newtheorem{lemma}  {\hspace*{1pc} Lemma}

\newtheorem{rem}  {\hspace*{1pc} Remark}
\newtheorem{prop}  {\hspace*{1pc} Property}

\begin{document}

\begin{frontmatter}



\title{A slot-based energy storage decision-making approach for optimal 
Off-Grid telecommunication operator}


\author[lab1]{Youssef Ait El Mahjoub}
\author[lab2,lab3]{Jean-Michel Fourneau}

\affiliation[lab1]{organization={Efrei Research Lab},
            addressline={Université Paris-Panthéon-Assas}, 
            city={Villejuif},
            postcode={94800}, 
            country={France}}

\affiliation[lab2]{organization={DAVID Laboratory},
            addressline={Université Paris-Saclay}, 
            city={Versailles},
            postcode={78000}, 
            country={France}}

\affiliation[lab3]{organization={Inria, ARGO},
            city={Paris},
            country={France}}

\begin{abstract}
This paper proposes a slot-based energy storage approach for decision-making in the context of an Off-Grid telecommunication operator. We consider network systems powered by solar panels, where harvest energy is stored in a battery that can also be sold when fully charged. To reflect real-world conditions, we account for non-stationary energy arrivals and service demands that depend on the time of day, as well as the failure states of PV panel. The network operator we model faces two conflicting objectives: maintaining the operation of its infrastructure and selling (or supplying to other networks) surplus energy from fully charged batteries. To address these challenges, we developed a slot-based Markov Decision Process (MDP) model that incorporates positive rewards for energy sales, as well as penalties for energy loss and battery depletion. This slot-based MDP follows a specific structure we have previously proven to be efficient in terms of computational performance and accuracy. From this model, we derive the optimal policy that balances these conflicting objectives and maximizes the average reward function. Additionally, we present results comparing different cities and months, which the operator can consider when deploying its infrastructure to maximize rewards based on location-specific energy availability and seasonal variations.
\end{abstract}



\begin{keyword}
Energy storage \sep Non Stationary Arrivals \sep Markovian Decision Process \sep  Empirical Data \sep Average Reward
\end{keyword}

\end{frontmatter}




\section{Introduction}
Despite the need for mobile infrastructures in developing countries, there are major challenges in providing power to these
systems, especially when they are Off-Grid. Thus energy harvesting and energy storage have to be deployed 
in the context of an off grid telecom infrastructure. With the recent increase in energy prices, the powering cost 
 has become the dominating operational cost for mobile networks operators
 \cite{PGMRD18}. Although the base stations for 5G and 6G mobile networks are designed to be energy efficient, their 
 large scale deployment will increase the demands of energy for the network infrastructure. 
 At the same time, deploying an energy harvesting system offers operators the opportunity to sell energy. Balancing two conflicting objectives, maintaining the reliability of their own network infrastructure and acting as an energy provider. This paper addresses the stochastic decision-making challenge raised by these conflicting objectives. 
 
We consider an energy storage system (i.e. a battery) which is filled by a sporadic random arrival of energy.
Following \cite{Gele12, Gele16, YHJM20} we represent energy by  discrete units called  Energy Packets (EP in this paper). 
This 
discretization of the energy allows to model the system as a discrete state stochastic system and leads in the literature to 
efficient analysis techniques based on closed form multiplicative solutions for the steady-state distribution 
of the system. 
An EP is typically the amount of energy needed to perform a job: typically sending a packet. 
The use of energy storage (mostly rechargeable batteries), even if it costly,   is justified by 
the erratic and unpredictable energy harvesting patterns \cite{DMRe16,ReMe19}. But these batteries also 
offer the opportunity for the operator to be sold to obtain a new source of  incomes and decrease the OpEx (Operating Expenses). 
The sporadic nature of the energy arrival requires some analysis of stochastic models to obtain the 
 distribution of the battery associated to energy harvesting devices (for instance solar panels)
and used to powered a mobile network node. 
In \cite{ALJM20}, the authors propose minute-scale models to estimate solar energy received by photovoltaic panels, based on solar irradiance and the clear sky index. In \cite{KGCC23}, a transient analysis of diffusion processes is conducted. The focus of this model was to calculate the probability that the battery depletes due to network operations and adverse solar conditions. In subsequent work \cite{KGCT23}, the same authors focuse on determining the probability density for the time required to charge batteries fully and the time to completely deplete stored energy. More recently, \cite{KGCT24} introduces a model for a solar energy harvesting and storage system that accounts for Day/Night duality in the context of green Off-Grid cellular base stations. The Day/Night duality has also been explored by some of us \cite{YAEM24} and in the work by Miozo \cite{MZDR14}, where the focus is on the performance evaluation of such systems under the variability inherent in stochastic processes affecting energy arrival and consumption. However, these studies do not address the decision problem associated with selling energy. 
 Here, we assume that when the battery contains a certain level of energy, it can be sold and it is 
 replaced by a new one which is empty. Thus the operator increases its income by selling a filled battery
 but it makes its infrastructure less reliable as the system may require energy before the solar panel can provide
 enough EP. To the best of our knowledge, as a continuation of our prior work \cite{YAEM24}, this paper is the first attempt to use the EP modeling assumption for a decision-making problem with a Markov Decision Process formulation (MDP in the following) with non-stationary arrivals.

In the development of the analysis, we found that the MDP exhibits a highly regular structure, which can be leveraged to simplify the numerical analysis of the optimal policy. Consequently, this paper makes three main contributions:
\begin{itemize}
    \item We address the decision-making problem faced by an operator who not only manages its infrastructure but also serves as a battery provider. By employing Markovian models, we solve this problem more efficiently than with traditional learning approaches, a result largely attributed to the independence assumption in energy production and service demands. 
\item We develop an enhanced algorithm for analyzing the MDP, allowing for faster processing and a larger range of possible actions and states. This algorithm’s efficiency stems from incorporating the graph structure of the MDP. 
\item We consider non-stationary arrivals for both energy production and service demands while considering PV panels failure states, aligning the model with real-world data and assumptions. This approach enables the telecommunications operator to make more informed decisions about infrastructure deployment across different locations and months of the year. This adaptation marks a significant advancement from our previous study in \cite{AYFJ24}, where arrivals were assumed to be stationary during the 'Day' period, without accounting for potential PV failure states.
\end{itemize}

The technical part of the paper is as follows. In Section II, we describe the model of a battery which receives energy from a sporadic source like a solar panel. This battery is used to feed a networking system but also the battery can be sold when it is sufficiently filled.  We assume that arrivals are random and described by an Interrupted Batch Process (IBP hereafter), with two phases corresponding to PV-ON and PV-OFF. PV-ON correspond to states where the PhotoVoltaic (PV) panel is operational while PV-OFF corresponds to failure states due to manufacturing defects or premature wear. Such a stochastic process is a discrete-time version of one of the stochastic processes proposed in \cite{MZDR14}. In this paper, the authors use solar measurements to build stochastic processes representing the energy arrival for photovoltaic sources and verify their accuracy in models of wireless systems.
We show that the Markov chain has a structure already defined in the literature \cite{Rob90, YHJM19} which is known to simplify the steady-state analysis of the chain. 
In Section III, we focus on defining, analyzing and solving the MDP and demonstrate how its structure can be used to find the optimal policy. To the best of our knowledge, this structure was only considered yet to 
analyze Markov chains. We define the reward system, which includes penalties for non-operational networking system and energy waste (i.e. energy packets lost), as well as positive rewards for the potential sale of batteries. Additionally, we introduce an algorithm that 
takes into account
the MDP structure. Finally, in Section IV  we study the numerical examples and we report the time needed to obtain the optimal policy with our method and with usual numerical techniques. We also provide results for different locations and months based on empirical open-access data.

\section{Battery Filling Markov Chain Model}
\subsection{The state}
We consider a discrete-time and discrete-space model where the state of the system is denoted by $(H, X, M)$. The component $H$ represents a clock ranging from $t_0$ to $T$, where $t_0$ marks the first hour of the day when the solar panel begins to receive energy, and $T$ represents the deadline for battery charging. In practice, $T$ can be interpreted as the last hour of the day when the photovoltaic (PV) panel receives energy.

The component $X$ represents the state of the battery, i.e., the amount of energy it contains, measured as an integer value of energy packets (EP). The battery's capacity is bounded between 0 and $C$, where $C$ is the maximum amount of energy the battery can store. This model is consistent with the Energy Packet Network, as introduced by Gelenbe and others \cite{Gele12, Gele16, GeZh19, YHJM20}, where the energy is discretized into EP units.

Finally, the component $M$ describes the operational status of the PV panels. It takes values in $\{\text{ON}, \text{OFF}\}$, where ON corresponds to the state where the solar panel is operational and supplying energy to the battery, and OFF represents the state where the panel is non-functional (i.e. a failure state), and thus no energy is supplied. The variable $M$ follows a two-state Markov chain, which modulates the transitions of the system. For example, when $M$ is in the OFF state, no solar energy is available to charge the battery. 

\subsection{The transitions}

\begin{figure}[!ht]
    \centering
    \includegraphics[width=13cm,height=7cm]{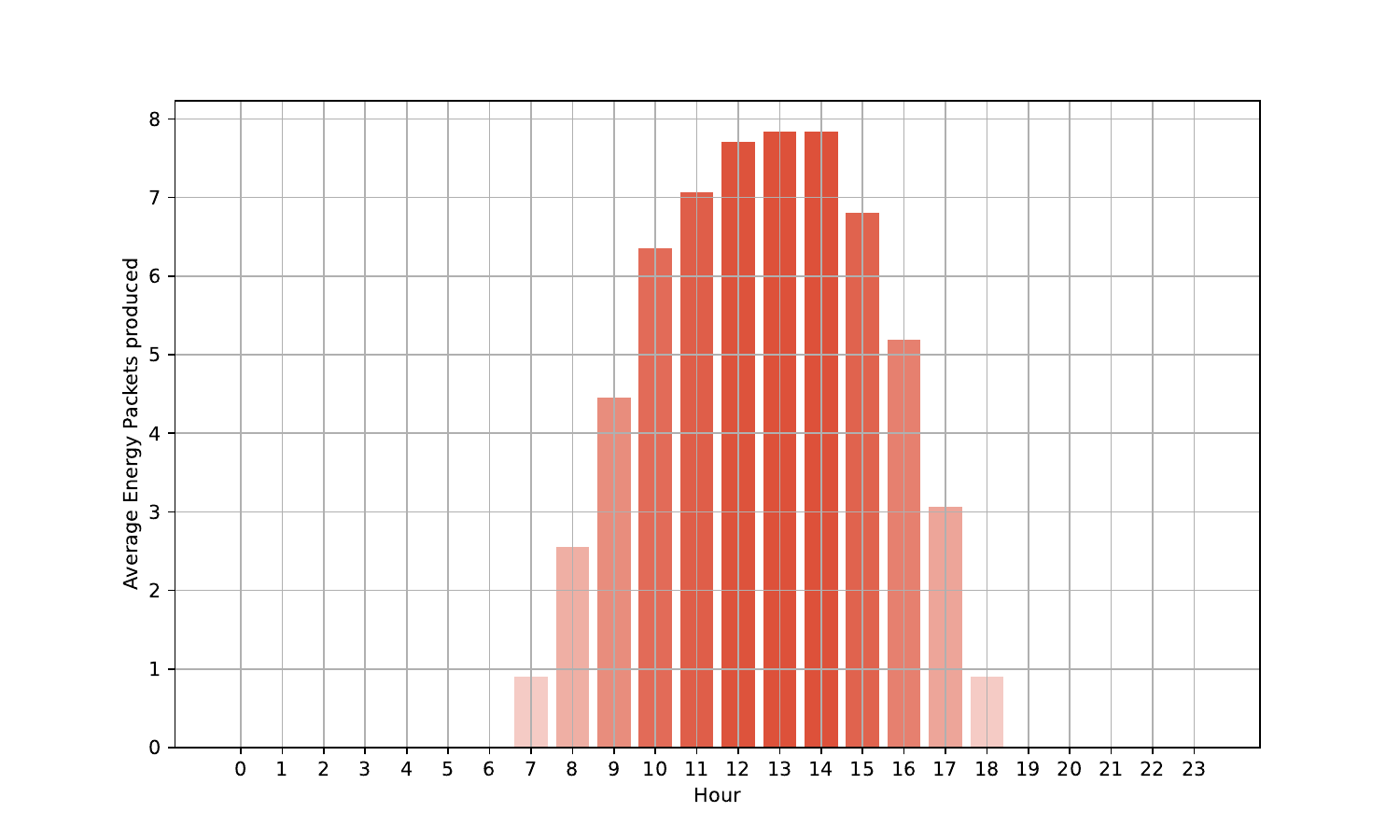}
    \caption{Average number of Energy Packets produced by a solar panel. Barcelona, August \cite{NREL}.}
    \label{EP-arrival}
\end{figure}

When the phase is equal to ON, the system evolves as follows. 
At each time step where the clock is strictly positive, the clock increases by one time slot. Then :

(i)  The battery receives a discrete amount of energy from the environment. The amount of energy is a batch with distribution $\mathcal{A}_H$, which depends on the hour hence the energy arrivals are non-stationary throughout the day. For example, the system may receive more energy during zenith hours and less towards the end of the day (see Fig. \ref{EP-arrival}).

\begin{figure}[!ht]
    \centering
    \includegraphics[width=12cm,height=6cm]{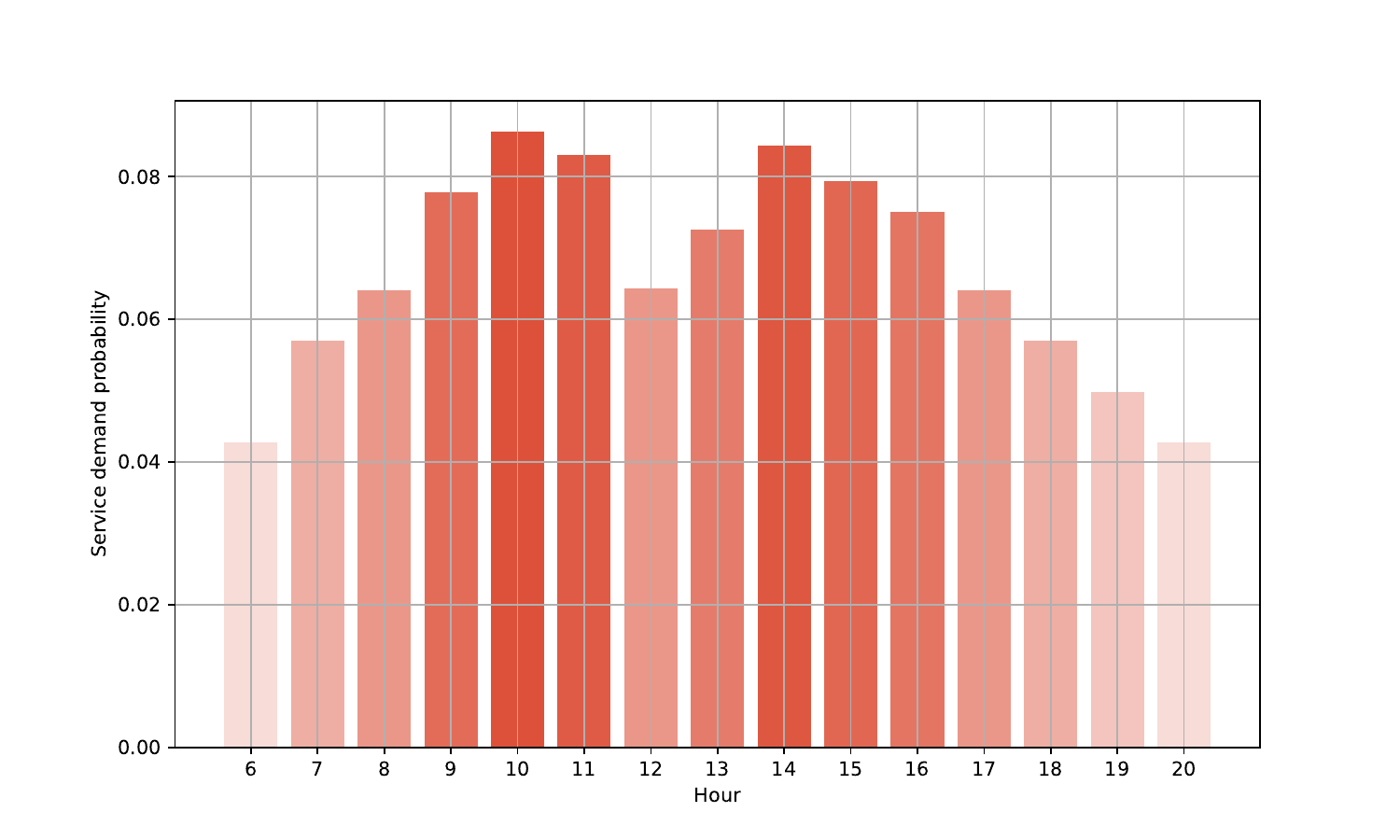}
    \caption{Distribution of Data Packets arrival i.e. service demands. Erlangs-based workload \cite{Cisco07}.}
    \label{DP-arrival}
\end{figure}

(ii) The operating system may require one EP (Energy Packet) to process a job, provided there is available energy in the battery. Otherwise, the job cannot be processed on time, and a penalty is incurred for this delay. The system consumes an EP according to a Bernoulli distribution ${\mathcal{B}_H}$, which also depends on the hour (see Erlang workload in Fig. \ref{DP-arrival}). Again, we assume non-stationarity and independence for data packets (DP) arrivals.


(iii) If the energy level in the battery exceeds its capacity $C$, any excess EPs from the incoming batch are lost. A penalty is applied for this loss of energy. Afterward, the battery may be released and sold to obtain a reward (details of which will be provided later).

We begin by specifying two conditions under which a battery can be released (represented by the green states in Figure \ref{fig1} (page 9):

\begin{itemize}
    \item The first condition is when the battery's energy level exceeds a threshold $F$. In this case, the battery is released with a probability that depends on the energy level but not on the clock. Let $d_x$ denote this probability when the system is in state $(h,x,\text{ON})$. This is modeled by a Bernoulli distribution ${\cal{Z}}_X^{ON}$, which takes the value $1$ when the battery is released. We assume independence.
    
    \item The second condition is when the clock reaches the deadline time $T$, irrespective of the battery's energy level.
\end{itemize}

Once the battery is released, the system jumps to the state $(0,0,\text{ON})$: the clock is reset, and a new empty battery is used. The release of a battery provides a reward that depends on the energy level. This reward is positive if the energy level is greater than or equal to $F$, and non-positive otherwise. Furthermore, we assume that the clock does not start increasing until the first energy arrival or a PV failure that could occur in the initial state. The time interval when the clock does not increase can be analyzed separately.

(iv) The PV panel may experience a technical issue, causing the system to enter in failure state PV-OFF. Thus, the PV panel can no longer produce energy. To simplify equation readability, we refer to PV-ON states as ON states and PV-OFF states as OFF states. The evolution between the two states of the modulating phase $M$ is rather simple: from ON to OFF 
the transition probability is $\alpha$ and the system stays at ON with probability $1-\alpha$. Similarly, the transition probability 
from OFF to ON is $\beta$ while the system stays at OFF with probability $1-\beta$. These transitions are associated with random variables distributed respectively according to ${\cal{Y}}_X^{OFF}$ and ${\cal{Y}}_X^{ON}$. The transitions of $M$ 
are represented by blue arcs in Figure \ref{fig1}  (page 9).

At each time slot (say $n$), depending on the state of phase $M_n$, 
we sample  some random variables. 
When $M_n = ON$, we have to consider $4$ independent  r.v. :  $e_n$ for EP arrivals, $b_n$ for DP service, $z^{ON}_n$ for battery release, and $\phi^{ON}_n$ for phase change
from distributions $\cal{A}_H$, ${\cal{B}}_H$,
${\cal{Z}}_X^{ON}$ and ${\cal{Y}}_X^{ON}$. We assumed that distributions  ${\cal{Z}}_X^{ON}$ and
${\cal{Y}}_X^{ON}$ depend on the value of the battery $X$.  
When $M_n = OFF$, we only have to consider $3$ independent  r.v.  $b_n$, $z^{OFF}_n$, and $\phi^{OFF}_n$
from distributions  ${\cal{B}}_H$, ${\cal{Z}}_X^{OFF}$ and ${\cal{Y}}_X^{OFF}$. We also assume that distributions ${\cal{Z}}_X^{OFF}$ and ${\cal{Y}}_X^{OFF}$ depend on energy level $X$.

The evolution of the system at time $n$  
when $H_n > 0$ and $M_n = ON$ is described by the following equations: \\ 
${\rm if~} (\phi^{ON}_n = 0)$ then $M_{n+1}  =  ON$, and $H_{n+1}$, $X_{n+1}$ are set to
 \begin{equation}
 \label{trans1}
\left[
\begin{array} {lr}
 {\rm if~} (H_n<T) \wedge  (X_n <F \vee z^{ON}_n=0)  & H_{n+1} =  H_n+1, \\
{\rm if~} (H_n=T)    \vee  (X_n  \geq  F \wedge z^{ON}_n=1)  & H_{n+1} = t_0, \\
{\rm if~} (H_n<T)  \wedge (X_n<F)  & X_{n+1} =  max(min(X_n + e_n,C)-b_n,0), \\
{\rm if~} (H_n=T) &  X_{n+1} =  0,  \\
{\rm if~} (H_n<T)  \wedge (X_n\geq F) \wedge (z^{ON}_n =1) & X_{n+1} =  0,  \\
 {\rm if~} (H_n<T)  \wedge (X_n\geq  F) \wedge (z^{ON}_n =0) & X_{n+1} =  max(min(X_n + e_n,C) -b_n,0). 
 \end{array}
\right. 
\end{equation}
${\rm if~} (\phi^{ON}_n = 1)$ then $M_{n+1}  =  OFF$, $H_{n+1} = H_n + 1$, $X_{n+1} = X_n$. \\
For instance, when $H<T$ and $X_n <F$, the probability of a transition from $(h,x,ON) $ to $(h+1,x-1,ON) $ is 
\[ {\cal{A}_H} [0]  {\cal{B}}_H [1] {\cal{Z}}_X^{ON} [0] {\cal{Y}}_X^{ON} [0]. \] As the random variables are independent, the resulting probability of this transition is the product of the probability that the phase does not change 
(i.e. with probability ${\cal{Y}}_X^{ON} [0] = (1-\alpha)$), with $0$ EP arrival probability (i.e. with probability ${\cal{A}_H} [0]$), and one DP service demand (i.e. with probability ${\cal{B}}_H [1]$), without releasing the battery (i.e. with probability ${\cal{Z}}_X^{ON} [0] $). 

The evolution of the chain when the system is in the particular starting state $(t_0,0,ON)$ is in the following. We assume that the size of the battery $C$
is larger than the largest arrival batch in distributions $\cal A_H$. This assumption allows to simplify the evolution equation 
to simplify the presentation but it is not needed for the numerical analysis.   
\\ - If  $e_n>0 $, then $H_n$ jumps to $t_0+1$ and $X_n$ jumps to $max(0,e_n-b_n)$. 
\\ - If $\phi_n = 1$, then $H_n$ jumps to $t_0 + 1$ and $M_n$ to OFF.
\\ - Otherwise we have a loop on state $(t_0,0,ON)$. 
\\ Note that at each time slot, both an arrival and service can occur. Hence an EP could be treated immediately after its arrival, which explains for instance the transition from $(t_0 = 9,0,ON)$ to $(10,0,ON)$ in Figure \ref{fig1}.

Let us now consider the case where $M_n=OFF$. There is no arrivals of EP but it is still possible to use the battery if is not empty. Thus the battery is always smaller than $C$ by induction: \\ 
${\rm If~} (\phi^{OFF}_n = 0)$ then $M_{n+1}  =  OFF$, and $H_{n+1}$, $X_{n+1}$ are defined as 

\begin{equation}
 \label{trans2}
\left[
\begin{array} {lr}
 {\rm if~} (H_n<T) \wedge  (X_n <F \vee z^{OFF}_n=0)  &  H_{n+1} =  H_n+1, \\
{\rm if~} (H_n=T)    \vee  (X_n  \geq  F \wedge z^{OFF}_n=1)  & H_{n+1} = t_0, \\
{\rm if~} (H_n<T)  \wedge (X_n<F)  & X_{n+1} =  max(X_n -b_n,0),  \\
{\rm if~} (H_n=T) &  X_{n+1} =  0,  \\
{\rm if~} (H_n<T)  \wedge (X_n\geq F) \wedge (z^{OFF}_n =1) & X_{n+1} =  0,  \\
 {\rm if~} (H_n<T)  \wedge (X_n\geq  F) \wedge (z^{OFF}_n =0) & X_{n+1} =  max(X_n -b_n,0).
 \end{array}
\right. 
\end{equation}
${\rm If~} (\phi^{OFF}_n = 1)$ then $M_{n+1}  =  ON$, $H_{n+1} = H_n + 1$ and $X_{n+1} = X_n$. \\
Finally, from state $(t_0, 0, \text{OFF})$, there is two possible transitions. Indeed, no EP are entering the system when the PV is in failure and it is not possible to serve a job as the battery is empty. Therefore we have a transition from $(t_0,0,\text{OFF})$ to $(t_0,0,\text{ON})$ with probability $\beta$ (dotted blue arc in Figure \ref{fig1}) and a loop at state $(t_0,0,OFF)$ with probability $1- \beta$. 
\subsection{The graph structure}
\begin{figure}[!ht]  
\hspace{-3.2cm}
\includegraphics[scale=0.5]{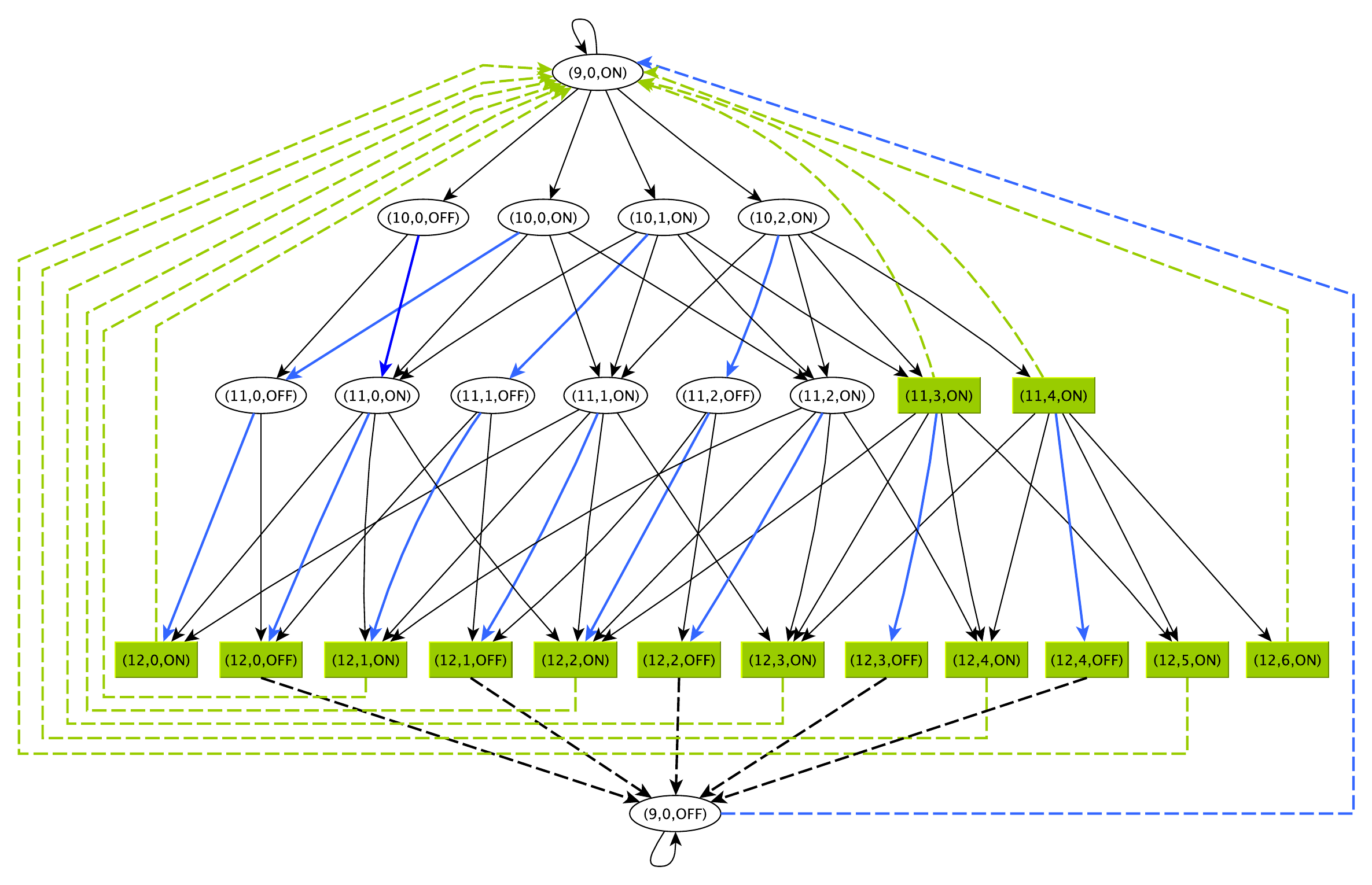}
        \caption{Toy example filling process with parameters $t_0=9h$, $T=12h$, $F=3$, $C=3$, and $\forall h \in \{t_0,\dots, T\}$, ${\cal{A}_H} = \{0, 1, 2\}$. \label{fig1}}
\end{figure}

\begin{prop} 
Based on the independence assumptions, the description of states and transitions $(H_n,X_n,M_n)_n$ is a Markov chain. 
\end{prop}
In Figure \ref{fig1}, we show a simple example of the Markov chain of this filling process.

In (\cite{Rob90}, page 141) Robertazzi has studied two Markov chains structures which can be used to solve efficiently the steady-state distribution of the chains based on these structures.  
Here we consider type B structure. 
\begin{defin} 

In a Robertazzi type B structure for a Markov chain, all the directed
cycles of the chain goes throught only one state. 
\end{defin} 

\begin{prop} 
Markov chain $(H_n,X_n,M_n)_n$ has a  Robertazzi  type B structure  where all the directed cycles of the Markov chain go 
through state $(t_0,0,ON)$.
\end{prop} 
\begin{proof}
The proof is based on the fact that the clock always 
increase when the battery is not empty. And the clock returns to $t_0$ when we release a battery.
More formally, emptying the battery leads to two states $(t_0,0,\text{ON})$ or $(t_0,0,\text{OFF})$ due to the modulating phase. 
When $M_n=ON$, the only state which have a self loop is state $(t_0,0,\text{ON})$ because when $H_n>0$ it increases at each transition of the Markov chain. Second, the same argument shows that one cannot have a directed cycle between state $(h1,x1,\text{ON})$
and $(h2,x2, \text{ON})$. Without loss of generality assume that $h1 \geq h2$.  As $H$ is increased at each transition, the directed path from $h1$ to $h2$ reset the clock and it goes through $(t_0,0,\text{ON})$.
When $M_n=OFF$, the property on the increasing clock during transitions (state $(t_0,0,\text{ON})$ and $(t_0,0,\text{OFF})$ excluded) is still valid. And from $(t_0,0,\text{OFF})$ there is only one transition that leads another state, the state $(t_0,0,\text{ON})$. 
Therefore all the cycles goes through state $(t_0,0,\text{ON})$. 
\end{proof}

Decision variables we will emphasize in next section are distributions ${\cal{Z}}^M_X$ ($M\in \{\text{ON},\text{OFF}\}$) and ${\cal{B}}_H$  for all values of capacity of battery level $X \geq F$. These distributions will depend on the action.

\section{Battery Filling Markovian Decision Process}
\begin{figure}[!ht]
    \centering
    \includegraphics[width=0.8\linewidth]{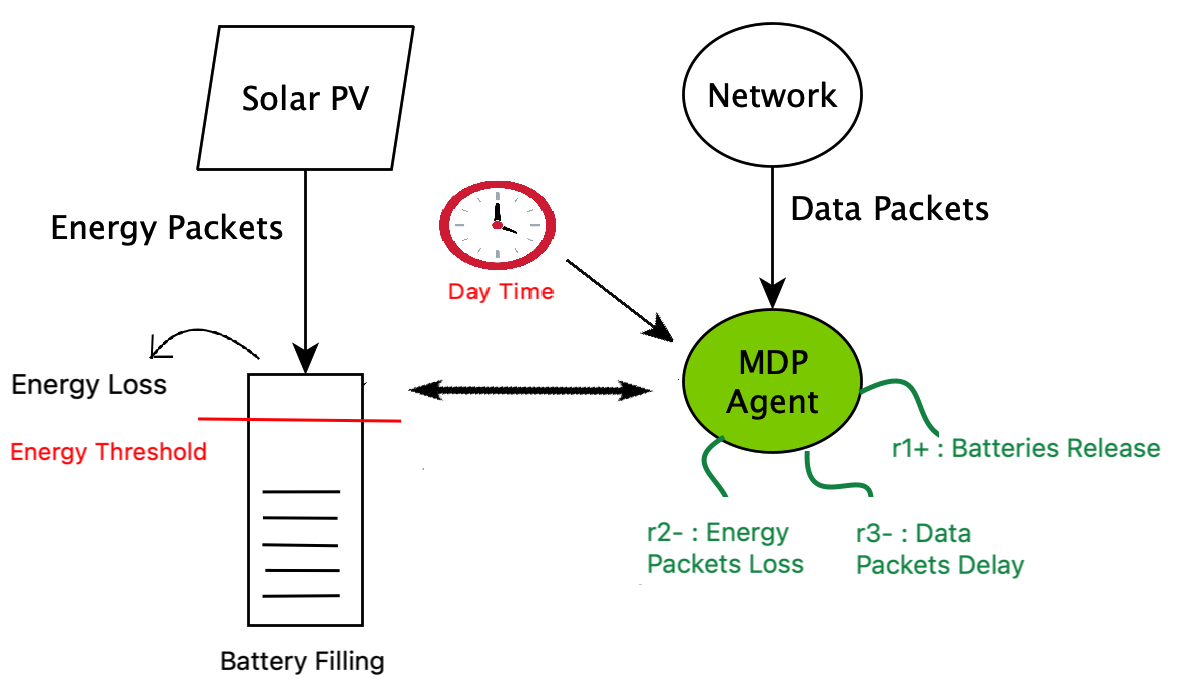}
    \caption{System representation}
    \label{figModel}
\end{figure}

With an Markovian Decision Process (MDP) formulation, we address the challenge of optimizing battery recharging to rapidly supply energy to the network while preventing the deployment (or sale) of depleted batteries. The system is equipped with an intelligent agent embedded within the battery. The main objective is to optimize the control parameters of this intelligent agent to ensure that battery supplies the network in an optimal manner. By strategically managing the recharge cycles and the energy output, the agent aims to maximize operational efficiency under constraints related to energy availability and time-sensitive demands. This MDP process aims to determine the optimal policy that guides the battery's behavior, balancing between maximizing battery release gains and minimizing the risk of data packets delays, while also considering potential energy packets losses. A schematic description is presented in Fig. \ref{figModel}.

\subsection{MDP formulation : states and actions }
Let $s=(h,x,m)$  represent a state of the battery model within the Discrete Time Markovian Decision Process $\{\mathbb{S},\ \mathbb{A}, \ P^{(a)},\ \mathbb{R}^{(a)} \}$. We define the state space as  $\mathbb{S} = \{h, x, m \in   \mathbb{N}  \ | \ t_0\leq h \leq T, \ 0\leq x \leq C,\ m\in\{ON,OFF\} \}$, and $\mathbb{A}$ a finite set of actions.
An action refers to the filling process of the previous section with similar EP arrivals distributions ${\cal{A}_H}$ in each action, but different service and release distributions. Let ${\cal{Z}}_X^{ON,(a)}$,  ${\cal{Z}}_X^{OFF,(a)}$ and ${\cal{B}}_H^{(a)}$  be, respectively, the Bernouilli distribution of a battery release in ON states, OFF states, and serving a task while considering action $a\in \mathbb{A}$. We assume that modulating phase probabilities $\alpha$ and $\beta$ do not depend on the action.
Consequently, the agent experiences non-stationary and identical arrivals for each action but exhibits different probabilities of either releasing the battery or immediately serving a job. The battery management strategy aims to control the optimal quantity of energy packets sent to the network, constrained by the need to avoid dispatching depleted batteries and observing time deadlines. 
\subsection{MDP formulation : probabilities per action}
In the following, we present transition matrix $P^{(a)}$ for each action. We first define the three set of states

\begin{equation}
 \left[
\begin{array}{l}
\mathbb{S}^m_1 = \{s=(h,x,m)\in \mathbb{S} \ |\  (0<h<T) \wedge (x<F) \}, \\
\mathbb{S}^m_2 = \{s=(h,x,m)\in \mathbb{S} \ |\  (0<h<T) \wedge (x\geq F) \}, \\
\mathbb{S}^m_3 = \{s=(h,x,m)\in \mathbb{S} \ |\  h = T \}.
\end{array}
\right.
\end{equation}

$\mathbb{S}^m_2$ and $\mathbb{S}^m_3$ refer to the releasing of battery states (i.e. green states in Fig. \ref{fig1}) in phase $m\in\{ON,OFF\}$. $\mathbb{S}^m_1$ refers to other states. For PV-ON states : \\
if $s\in \mathbb{S}^{ON}_1$ then
\begin{equation}
\label{equP1}
 \left[
\begin{array}{lr}
    if~ s'=(h+1,max(min(x + e_n,C)-b_n,0), \text{ON}) &  P_{s,s'}^{(a)} =  (1-\alpha){\cal{A}_H} [e_n]  {\cal{B}}_H^{(a)} [b_n],  \\
     if~ s'=(h+1,x,OFF) &  P_{s,s'}^{(a)} = \alpha.
\end{array}
\right.
\end{equation}
if $s\in \mathbb{S}^{ON}_2$ then
\begin{equation}
\label{equP2}
 \left[
\begin{array}{lr}
    if~ s'=(h+1,max(min(x + e_n,C)-b_n,0), ON) & P_{s,s'}^{(a)} =  (1-\alpha){\cal{A}_H} [e_n]  {\cal{B}}_H^{(a)} [b_n] {\cal{Z}}_X^{ON,(a)} [0],   \\
     if~ s'=(h+1,x,OFF) &   P_{s,s'}^{(a)} = \alpha, \\
     if~ s'=(t_0, 0, ON) &   P_{s,s'}^{(a)} = (1-\alpha){\cal{Z}}_X^{ON,(a)}[1].
\end{array}
\right.
\end{equation}
and if $s\in \mathbb{S}^{ON}_3$ we have $P_{s,s'}^{(a)} = 1$. \\ It remains to address the special case $(t_0,0,ON)$, as the process of filling begins when the first non-zero arrival occurs or when PV panel exhibits a failure event. Hence the chain remains at $(t_0,0,ON)$ with probability $(1-\alpha){\cal{A}_H} [e_n=0]$ and other transitions have probabilities: $(1-\alpha){\cal{A}_H} [e_n>0]  {\cal{B}}_H^{(a)}[b_n]$ for transitions in PV-ON states, while having $\alpha$ as probability for $(t_0,0,OFF)$ state transition.

The same reasoning is applied to obtain PV-OFF states probabilities : \\
if $s\in \mathbb{S}^{OFF}_1$ then 
\begin{equation}
\label{equP3}
 \left[
\begin{array}{lr}
    if~ s'=(h+1,max(x -b_n,0), OFF) &  P_{s,s'}^{(a)} =  (1-\beta) {\cal{B}}_H^{(a)} [b_n],  \\
     if~ s'=(h+1,x,ON) &  P_{s,s'}^{(a)} = \beta.
\end{array}
\right.
\end{equation}
if $s\in \mathbb{S}^{OFF}_2$ then
\begin{equation}
\label{equP4}
 \left[
\begin{array}{lr}
    if~ s'=(h+1,max(x-b_n,0), OFF) &  P_{s,s'}^{(a)} =  (1-\beta) {\cal{B}}_H^{(a)} [b_n] {\cal{Z}}_X^{OFF,(a)} [0],  \\
     if~ s'=(h+1,x,ON) &  P_{s,s'}^{(a)} = \beta, \\
     if~ s'=(t_0, 0, OFF) &  P_{s,s'}^{(a)} = (1-\beta){\cal{Z}}_X^{OFF,(a)}[1].
\end{array}
\right.
\end{equation}
and if $s\in \mathbb{S}^{OFF}_3$ we have $P_{s,s'}^{(a)} = 1$. \\ 
It remains to address the special state $(t_0, 0, OFF)$. This one clearly goes to $(t_0,0,ON)$ with probability $\beta$, when the deadline has been achieved. Otherwise, the chain stays at $(t_0,0,OFF)$ with probability  $1-\beta$. One can verify that for all states (notably,  Equation \eqref{equP1}, \eqref{equP2}, \eqref{equP3}, and \eqref{equP4}) the sum of transition probabilities is equal to $1$.
\subsection{MDP formulation : instant rewards}
We represent immediate rewards formula $r(s,a)$, representing rewards the system receives upon executing action $a$ from state $s$ as 
\begin{equation}
    r(s,a)= \sum_{s'}P_{s,s'}^{(a)}.\mathbb{R}^{(a)}_{s,s'}
\end{equation}
and 
\begin{equation}
\label{equaReward}
\mathbb{R}^{(a)}_{s,s'} =  \left[
\begin{array}{ll}
g(x).r^+_1 & if \ (h' = 0)  \\
max(0, x+e_n-b_n-C).r^-_2 & if \ (x' = C)  \\
r^-_3 & if \ (x' = 0)  
\end{array}
\right.
\end{equation}
Where $r^+_1$ represents a positive reward that the agent receives for each energy packet sold at the moment of battery release (i.e. when $h'=0$). $g(x)$ is a function assumed to be linear, which thus promotes releasing not only above a threshold but also close to the battery capacity. $r^-_2$ denotes a negative reward that the agent receives when a packet is lost, meaning that the current battery cannot accommodate it (i.e. when $x'=C$). Given that arrivals occur in batches, some packets from the batch can be stored if capacity allows, but the remaining packets are counted as lost. This accounts for the multiplicative term with  $r^-_2$. Finally, we consider $r^-_3$ penalty, which is applied when a transition results in an empty battery (i.e., when $x'=0$). This condition means that the operating system related to the battery cannot serve a task, hence the data packet is delayed. These three rewards collectively reflect the trade-off between the need to sell batteries and the risks associated with leaving the battery empty, thereby preventing job services for the related operational network, and the risk of overfilling, which could lead to penalties due to energy packets loss.

\subsection{MDP resolution : average reward criteria}
We aim to identify the optimal policy for the agent to maximize the reward function across a prolonged sequence of decisions. Specifically, our focus is on maximizing the average reward over that trajectory.
\begin{prop} 
\label{prop_MDP_Rob}
For any policy $\pi$, the induced graph of $P^{(\pi)}$, of the defined DTMDP (Discrete Time Markovian Decision Process), has a Robertazzi type B structure.
\end{prop}
\begin{proof}
    The decision strategy of our model is to find the optimal probability of either releasing the battery above a threshold $x\geq F$ (i.e. with probability ${\cal{Z}}_X^{m,(a)}>0$), or  serving a task (i.e.  with probability ${\cal{B}}_H^{(a)}>0$). Hence transitions defined in Equation \eqref{trans1} and \eqref{trans2} (and the ones for particular states remains the same for any action taken by the agent. Only their probability changes with the action $(a)$ as stated in Equations \eqref{equP1}, \eqref{equP2}, \eqref{equP3}, \eqref{equP4}. 
\end{proof}

In following, we leverage this structural pattern to evaluate efficiently, with exact results, any policy in the relative evaluation step of the relative policy iteration algorithm in average reward MDPs \cite{Putr94,Abhj15}. 

First, we recall Bellman equations in the context of Policy Evaluation algorithm. 
To optimize an average reward criteria, under policy $\pi$, one need to estimate a value function for all states and the average reward $\rho^{(\pi)}(s)$. Let $r(s_t,\pi(s_t))$ be the immediate reward obtained from state $s$ at time "t" taking action $\pi(s_t)$. Then 

\begin{equation}
\rho^{(\pi)}(s) = \lim_{T \to \infty} \frac{\sum_{t=1}^T r(s_t,\pi(s_t)) \ /s_t = s}{T}.
\end{equation}
For the average reward criteria, defined in the last equation, a simplification could be made if the MDP is unichain  \cite{Putr94}. That is the average reward does not depend on the state. In unchain policies, the induced graph generates a single recurrent class (with some transient states). Hence, states will be revisited indefinitely which leads, asymptotically, to similar average reward. Unlike multi-chain policies, which can generate multiple recurrence classes, resulting in a possible distinct average value for each recurrence class.
\begin{lemma}
\label{lemUnichain}
    The DTMDP is unichain. Therefore 
    \begin{equation}
    \forall s, s', \ \ \rho^{(\pi)}(s) = \rho^{(\pi)}(s') = \rho^{(\pi)}.
\end{equation}
\end{lemma}
\begin{proof}
   We first recall the following assumptions about distributions: probabilities of arrivals, service, and battery release, along with parameters $\alpha$ and $\beta$, are all non-negative but strictly less than $1$.
    Based on these assumptions and the system’s evolution equations, the DTMC associated with each action is irreducible (every state can be reached from every other state). The DTMC is therefore ergodic as it is irreducible and includes at least one loop (for instance, the state $(t_0,0,ON)$) which ensures it is aperiodic. Hence, given the ergodicity of the DTMC for each action, we can state that under any stationary policy $\pi$, all states belong to a single closed communicating class. Thus, the defined DTMDP is unichain. Therefore, the average reward obtained from a decision trajectory starting from any state within this closed communicating class converges to the same average value, $\rho^{(\pi)}$.
\end{proof}
Another concern is to define an estimation of each state s given the evolution model. This estimation is refered to as a value function $V^{(\pi)}$. However, the natural value function in average reward criteria tends to have difficulties with large values (contrary to discounted reward, where discount factor $\gamma<1$ assures to have bounded values from estimated future rewards). Hence a natural version of the value function is defined as, $\forall s\in S$  
\begin{equation}
     V^{(\pi)}(s) = \lim_{T \to \infty} \mathbb{E}^{(\pi)} \big[ \sum_{t=1}^T r(s_t,\pi(s_t)) \\ /s_t = s \big].
\end{equation}
The relative value function consists in retrieving the value function of some defined state "x" (i.e. the bias value or relative value), solving the magnitude problem \cite{Abhj15},  $\forall s\in S$

\begin{equation}
     V^{(\pi)}(s) = \lim_{T \to \infty} \mathbb{E}^{(\pi)} \big[ \sum_{t=1}^T r(s_t,\pi(s_t)) \ /s_t = s \big] -  V^{(\pi)}(x)
\end{equation}
\begin{equation}
\label{eq1}
      \Rightarrow \ \ V^{(\pi)}(s) = r(s,\pi(s)) -  \rho^{(\pi)} +  \sum_{s'=1}^{|\mathbb{S}|}P^{(\pi)}_{s,s'} \ V^{(\pi)}(s').
\end{equation}

This last equation is the Bellman equation for relative  policy evaluation \cite{Putr94} which consists on a system of $|\mathbb{S}|$ linear equations. The unknowns are vector $V^{(\pi)}$ and  scalar $\rho^{(\pi)}$. That could be either solved by classical linear solvers that comes with significant computational cost  or iteratively, with some lack of precision, using fixed point methods. One note that if the steady-state distribution for some policy, we note $\Pi^{(\pi)}$, exists then we can derive the average Markov reward process formula
\begin{equation}
\label{eq2}
    \rho^{(\pi)} = \sum_{s\in \mathbb{S}} \Pi^{(\pi)}(s).r(s,\pi(s)) .
\end{equation}
Once a policy is evaluated, one can use following equations to improve the policy. These equations use the value vector obtained from Equation \eqref{eq1}. The Q-function is defined as 
\begin{equation}
\label{eq3}
Q(s,a) = r(s,a) + \lambda \sum_{s'=1}^{|\mathbb{S}|}P^{(a)}_{s,s'}V(s'),
\end{equation}
hence optimal policy \cite{Putr94,Abhj15} in each state is defined as
\begin{equation}
\label{eq4}
    \pi(s) \in \arg\max_{a\in A(s)}\Big[ Q(s,a) \Big].
\end{equation}
The Relative Policy Iteration Algorithm, then, consists in starting with an arbitrary policy, evaluate it using Equation \eqref{eq1} and improve it, if possible, with Equation \eqref{eq3} and \eqref{eq4}. The algorithm stops once no improvements are possible in Equation \eqref{eq4}. Therefore that last policy is the optimal policy $\pi^*$. 
Our concern is to use the graph structural property to solve efficiently Equation \eqref{eq1}. First, we order the states in ascending order of $H$ (as shown in Figure \ref{fig1}). Hence, we denote the root state $(t_0,0,ON)$ as $s_1$, other states as $s_k$ with k$>$1. Obtained transition matrix $P^{(\pi)}$ induced by policy $\pi$  has the following structure
$P^{(\pi)} = C^{(\pi)} + U^{(\pi)} $ where $C^{(\pi)}$ is a matrix whose first column is positive and all the other entries are set to $0$, and $U^{(\pi)}$ is a strictly upper diagonal matrix.

Let $V^{(\pi)}$, $R^{(\pi)}$ and $e$ be row vectors of, respectively, values, immediate rewards and vector where all element set to $1$. Then, Equation \eqref{eq1} is expressed as 
\begin{equation}
    V^{(\pi)} = R^{(\pi)} - \rho^{(\pi)} e + V^{(\pi)} . P^{(\pi)^t}
\end{equation}
\begin{equation}
\label{eqV}
   \Rightarrow V^{(\pi)}.(I - (C^{\pi}   + U^{\pi})^t)  = R^{(\pi)} - \rho^{(\pi)} e.
\end{equation}
Where $P^{(\pi)^t}$ is the transposed matrix of $P^{(\pi)}$ and "." is the vector matrix product.

From graph structure, we derive the following lemma.
\begin{lemma}
\label{lemmaV}
    Let $\mathbb{C}$ be the set of states having non zero values in first column i.e.  $\mathbb{C} = \{ s\in \mathbb{S}, \  C^{(\pi)}[s,s_1] = 1 \}$. Then  
\begin{equation}
\label{eq7}
    V^{(\pi)}(s_1) = 0, 
\end{equation}
\begin{equation}
\label{eq8}
   V^{(\pi)}(s) = R^{(\pi)}(s) - \rho^{(\pi)}, \ \ \ \ \ \ \   \forall s\in \mathbb{C} 
\end{equation}
\begin{equation}
\label{eq9}
  V^{(\pi)}(s) = R^{(\pi)}(s) - \rho^{(\pi)} + \sum_{s'>s}^{|\mathbb{S}|} V^{(\pi)}(s') U^{(\pi)}[s,s'], \ \    \forall s\in \mathbb{S}  \smallsetminus \mathbb{C} .
\end{equation}
\end{lemma}
\begin{proof}
First, Property \ref{prop_MDP_Rob} asserts that for any policy $\pi$, the induced graph structure conforms to a Robertazzi type B structure. Second, Lemma \ref{lemUnichain} establishes that the Markov Decision Process is unichain. We then prove, how to simplify the resolution of Equation \eqref{eqV}:
\begin{itemize}
\item We assume that the bias state in relative policy evaluation phase is state $s_1$. Consequently, its value will be retrieved from the estimation values of the whole vector $V^{(\pi)}$. Thus $V^{(\pi)}(s_1) = 0$. 
\item From the construction of the graph, we define $\mathbb{C} = \{ s\in \mathbb{S}, \  C^{(\pi)}[s,s_1] = 1 \}$. The particular states in $\mathbb{C}$ have only one outgoing arc with probability $1$. Hence, for all $s\in \mathbb{C}$, Equation \eqref{eq1} can be simplified to
\begin{equation}
    V^{(\pi)}(s) = V^{(\pi)}(s_1) + R^{(\pi)}(s) - \rho^{(\pi)} e.
\end{equation}
Given that bias state is $s_1$ and its value is zero, the above equation simplifies to Equation \eqref{eq8}.
 \item For other states $s\in \mathbb{S}\smallsetminus \mathbb{C}$, Equation \eqref{eq9} is derived from the upper diagonal structure of the transition matrix of the graph. This structural property justifies the summation that begins with terms where $s'>s$. 
\end{itemize}
\end{proof}
From this construction, one need to develop, in sequence, Equation \eqref{eq7}, followed by Equation \eqref{eq8}, then equation \eqref{eq9}. In \eqref{eq9} one need to evaluate the states in descending order (i.e. from bottom states in the graph to the root state $s_1$).

\begin{rem}
In Lemma \ref{lemmaV}, we showed how to solve efficiently Equation \eqref{eq1} and obtain vector $ V^{(\pi)}$. However, $\rho^{(\pi)}$  is still unknown in Equation \eqref{eq8}. Here, we will use Markov reward process formula of Equation \eqref{eq2}, to derive the average reward $\rho^{(\pi)}$ for a policy $\pi$ which needs the steady state distribution probability.  Calculating steady-state probabilities can be particularly challenging when scaling: often balancing between computational accuracy over significant execution time or fast computations with lack of precision. One may adopt GTH (Grassmann, Taksar, and Heyman) algorithm, Power method, Power method + Gauss-Seidel iterations, among others \cite{Stew94}. 
 However for this graph structure, we showed in \cite{YHJM18,YHJM19}, that we can compute exact results efficiently in $O(m)$ where $m$ is the number of arcs in the graph. In the following Lemma, we propose a similar approach to compute that steady-state probability distribution.
\end{rem}
\begin{lemma}
\label{lemmaPi}
For all $s>s_1$, if $P^{(\pi)}$ is ergodic, then 
\begin{equation} 
\label{eq11} 
\Pi^{(\pi)} (s)  = \alpha (s) \ \Pi^{(\pi)} (s_1)
\end{equation} 
where $\alpha (s) $ is defined as
\begin{equation} 
\label{eq12}
 \left\{ 
    \begin{array}{l} 
       \alpha (s_1) = 1,  \\
       \alpha (s)  = \sum_{s'<s} \alpha(s') U^{(\pi)} [s',s],  \ \ \forall s>s_1.
    \end{array}
    \right.
\end{equation}
Furthermore, from normalisation of probabilities, we have
$\Pi^{(\pi)} (s_1)  = \big[ 1+ \sum_{s>s_1}^{|\mathbb{S}|}\alpha (s) \Big]^{-1} $.
\end{lemma}
\begin{proof}
The proof is inspired by our previous work \cite{YHJM19}. It derives immediately from the balance equations.
First, one need to suppose the ergodicity condition for transition graph of policy $\pi$. Which guarantees existency of the steady state probability distribution $\Pi^{(\pi)}$. Hence we can state the classical formulation 
\begin{equation}
\label{eq14}
\left[
\begin{array} {l}
\Pi^{(\pi)} P^{(\pi)} = \Pi^{(\pi)} \\
\Pi^{(\pi)} e^t = 1
\end{array}
\right. 
\Rightarrow
\left[
\begin{array} {l}
\Pi^{(\pi)} (C^{(\pi)}+ U^{(\pi)}) = \Pi^{(\pi)} \\
\Pi^{(\pi)} e^t = 1
\end{array}
\right. 
\end{equation}
that provides balance equations in each state. By induction on $s>1$, we have :
\begin{itemize}
\item For $s=2$, from Equation \eqref{eq14}, we clearly have:
$\Pi^{(\pi)} (2) = \Pi^{(\pi)}(1) P [1,2] $, 
because node $1$ is the only predecessor of node $2$. We clearly get the value for $ \alpha(2)$.
\item For an arbitrary $s>2$, assume the induction holds until $s-1$. We have: 
\[ 
\Pi^{(\pi)} (s) = \sum_{s'<s} \Pi^{(\pi)} [s'] \ P^{(\pi)} [s',s] = \sum_{s'<s} \Pi^{(\pi)} [s']\  U^{(\pi)} [s',s] . 
\] 
Using the induction for  $ \Pi^{(\pi)} (s') $ with $s'<s$, we get: 
\[ 
\Pi^{(\pi)} (s) = \sum_{s' <s}  \big[ \alpha(s')  \Pi^{(\pi)}(1) \ \big] \ U^{(\pi)} [s',s],
\] 
from which one can readily obtained the induction (as $\Pi^{(\pi)}(1)$ is initialized to $1$). 
\item Finally, from $\Pi^{(\pi)} e^t = 1$,  we compute the real value of 
$\Pi^{(\pi)} (1) $ from the expression:
\[
\Pi^{(\pi)} (1) +   \sum_{s>1} \Pi^{(\pi)}(s) = \Pi^{(\pi)} (1) \ \big[1 +   \sum_{s>1} \alpha(s) \ \big] = 1
\]
\end{itemize}
\end{proof}
\begin{rem}
The proposed evaluation phase methodology is not limited to our particular model but applies to any unichain MDP formulation that adheres to structure in Property \ref{prop_MDP_Rob}. The computational complexity of the evaluation phase is upper bounded by $O(2m)$. This upper bound comprises $O(m)$ for Lemma \ref{lemmaPi} and an additional $O(m)$ for Lemma \ref{lemmaV}. Therefore, it is highly suitable for sparse matrices.
\end{rem}

\subsection{MDP resolution : optimal policy measures }
From optimal policy $(\pi^*)$ obtained from precedent sections, we can derive the following measures.
Let $m\in\{\text{ON,OFF}\}$ and let $\Pi^{(\pi^*)}$ be the steady state probability distribution obtained from the graph induced by policy $(\pi^*)$. This distribution can be obtained from Lemma \ref{lemmaPi}.

\begin{itemize} 
\item The average number of EP stored in the battery at the release (say $\mathbb{E}[Release]$), which is decomposed into two terms: the first one corresponds to the reward
when the clock is equal to deadline $T$
and the second one is associated with 
the reward obtained by releasing a battery with a capacity larger than the threshold $F$. Battery gain may depend
on the capacity (it is a function of $x$, say $g(x)$). $g(x)$ may be non positive if $x < F$. Then $\mathbb{E}[Release] =$ 
\begin{equation}
\label{Erelease}
\sum_x \sum_m  \Pi^{(\pi^*)}(T,x,m) g(x)
+ 
\sum_{h<T} \sum_{x \geq F} \sum_m \Pi^{(\pi^*)}(h,x,m) g(x) {\cal{Z}}_x^{m,(\pi^*)}[1]. 
\end{equation}
\item Delayed data packets probability, i.e. when the battery is empty while having a service demand.
\begin{equation}
\label{EDelay}
\mathbb{E}[Delay] = \sum_{h} \sum_{m}  \Pi^{(\pi^*)}(h,0,m){\cal{B}}_h^{(\pi^*)} [1].
\end{equation}

\item Average number of lost EP when the battery is full. This measure is only triggered by PV-ON states that are subject to EP arrivals. Let $m_h$ be the maximal size of the batch of EP arrivals at time $h$, then $\mathbb{E}[Lost]=$
\begin{equation}
\label{ELost}
\sum_h \sum_{x \geq m_h} \sum_{e_n=0}^{m_h} \sum_{b_n=0}^1 \Pi^{(\pi^*)}(h,x,ON){\cal{A}_H} [e_n] {\cal{B}_H} [b_n] max(0, x+e_n-b_n-C).
\end{equation}
\end{itemize} 

\begin{algorithm}[p]
\SetAlgoLined
\SetKwInOut{Input}{Input}
\SetKwInOut{Output}{Output}
 \caption{\label{algo1}Steady-state probability algorithm for policy $\pi$}
 \Input{Transition matrix $P^{(\pi)}$}
 \Output{Steady-state probability distribution vector $\Pi^{(\pi)}$}
 \BlankLine
     Initialize $\alpha(s_1) = 1$   \\
     Using Equation \eqref{eq12}, get values of $\alpha(s)$ for all $s>s_1$ \\
     Deduce the value of $\Pi^{(\pi)}(s_1)$ from normalisation  \\
     Obtain $\Pi^{(\pi)}(s)$ for all $s>s_1$ using Equation \eqref{eq11}.
\end{algorithm}
\begin{algorithm}[p]
\SetAlgoLined
\SetKwInOut{Input}{Input}
\SetKwInOut{Output}{Output}
 \caption{\label{algo2}Modified Policy Iteration algorithm}
 \Input{State space $\mathbb{S}$; action space $\mathbb{A}$; transition matrices $P^{(a)}$, reward matrices $\mathbb{R}^{(a)}$.}
 \Output{Optimal policy $\pi^*$, value function $V^{(\pi^*)}$, combined average rewards $\rho^{(\pi^*)}$, measures $\mathbb{E}[Release]$, $\mathbb{E}[Delay]$, $\mathbb{E}[Lost]$.}
 \BlankLine
    Set $k \gets 1$   \\
    Select an arbitrary policy $\pi_k$ \\
    \textbf{Policy evaluation} : \algorithmiccomment{Lemma \ref{lemmaV}} \\
    - Use Algorithm \ref{algo1}, to obtain $\Pi^{(\pi)}$ 
    \algorithmiccomment{Lemma \ref{lemmaPi}} \\
    - Calculate $\rho^{(\pi)}$ using Equation \eqref{eq2} \\
    - Calculate $V^{(\pi)}_k(s_1)$ and $V^{(\pi)}_k(s)$ where $s \in  \mathbb{C}$, using Equation \eqref{eq7} and \eqref{eq8} \\
     \For{$s \in \mathbb{S} \smallsetminus \mathbb{C}$, by descending order }{
        - Calculate $V^{(\pi)}_k(s)$ from Equation \eqref{eq9} \\
     }
     \textbf{Policy Improvement} :  \\
     - Compute the Q-value from Equation \eqref{eq3} \\
     - Choose a new policy $\pi_{k+1}$ using Equation \eqref{eq4} \\
     \textbf{Stopping criteria} :  \\
     \eIf{$\pi_{k+1}(s) = \pi_{k}(s)$, $\forall s \in \mathbb{S}$}{
     -Set $\pi^*(s) \gets \pi_{k}(s)$ \\
     -Calculate $\mathbb{E}[Release]$, $\mathbb{E}[Delay]$, $\mathbb{E}[Lost]$   from \eqref{Erelease}, \eqref{EDelay} and \eqref{ELost}. \\
     -The algorithm stops.
     }
     {
     Set $k \gets k+1$, and go to \textbf{Policy evaluation} step.
     }
\end{algorithm}
\newpage
\section{Numerical results}
\subsection{Algorithms comparison}
In the following, we compare execution time of the proposed modified Relative Policy Iteration Algorithm (say RPI+RB). This algorithm primarily enhances the speed and accuracy of the policy evaluation phase. The policy evaluation, referred to by Equation \eqref{eq1}, can be either solved  iteratively using fixed point method (say RPI+FP) or using a direct method as Gauss-Jordan elimination (say RPI+GJ). We also compare results with the classical Relative Value Iteration (RVI). For iterative algorithms (RVI and RVI+FP), stopping criteria are either when $span(V^{\pi}_k - V^{\pi}_{k+1})<\epsilon$ or when a maximum number of iterations is reached. To test scalability, we present results from three scenarios ranging from small to large scale MDPs. In the first scenario (Table \ref{tabA}), we considered $10$ actions with the state space expanding to $5.10^2$. In the second scenario (Table \ref{tabB}), we increased the number of actions to $50$, with the state space reaching up to $8.10^3$. The final scenario (Table \ref{tabC}) includes $100$ actions, where each action covers up to $2.10^5$ states. We fixed $\epsilon = 10^{-10}$ and $MaxIteration=10^5$. From these results we observe that: 
\begin{itemize}
    \item The Relative Value Iteration RVI method yields acceptable results but suffers from slow convergence, necessitating further iterations compared to policy iteration algorithms. This is enhanced by the highly sparse structure of our graph, which is unsuitable for this algorithm. Results are notably contingent on the value of $\epsilon$, potentially failing to converge to the optimal solution if $\epsilon$ is not sufficiently small, yet also risking non-termination if $\epsilon$ is exceedingly small.
    \item The evaluation phase in RPI+GJ is performed by a Gauss-Jordan elimination which has a cubic complexity. That explains difficulties to handle large scale scenario.
    \item The proposed algorithm RPI+RB demonstrates the fastest results among the investigated algorithms. The evaluation phase performs a direct method and exhibits high numerical accuracy. In contrast with evaluation phase in RPI+FP which is iterative and relies on convergence criteria  $\epsilon$. For instance, in the case of the largest MDP we investigated, each action among $|\mathbb{A}|=100$ generates a model with $|\mathbb{S}|=2.10^5$ states and approximately $m\approx 8.10^5$ arcs. This MDP was solved with accurate and exact results in $3528$ seconds, compared to $6158$ seconds (in RVI algorithm) for an approximate solution with an accuracy of $10^{-10}$. Execution time of other algorithms is clearly inefficient, exceeding $10^4$ seconds.
\end{itemize}
Note that MDP models are generated using  Xborne tool (a tool that we have previously developed \cite{XBN16}) that can produce very large-scale sparse probabilistic models in $C$ language. Following this, we solve MDP models using a Python-based framework we developed, which efficiently handles sparse matrices through vectorizations. Computations have been performed on a laptop with $10$ cores ($8$ of them at $3.2$ GHz peak frequency and two others at $2$ GHz peak frequency) with $16$GB RAM.
\vspace{-0.2cm}
\begin{table}[!ht]
 \centering
   \caption{Execution time (in seconds) and number of iterations.\\ Small scale MDPs with $|\mathbb{A}| = 10$ actions. \label{tabA}}
\begin{tabular}{p{1.9cm}||c|c|c|c}
  & $|\mathbb{S}|=10^2$ & $|\mathbb{S}|=2.10^2$ & $|\mathbb{S}| = 3.10^2$ & $|\mathbb{S}|=5.10^2$   \\\hline \hline
RVI        & $0.60 \ | \ 10348$ & $0.68 \ | \ 9060$ & $0.84 \ | \ 9389$ &  $1.30 \ | \ 10148 $  \\\hline \hline
 RPI + FP  & $15.11 \ |\ 4$     & $19.60 \ |\ 3 $    & $28.18 \ |\ 3$    &  $88.72 \ |\ 5$ \\\hline
 RPI + GJ  & $0.12 \ |\ 4$      & $0.32 \ |\ 3$      & $0.64 \ |\ 3$     &  $3.27 \ |\ 5$ \\\hline
 RPI + RB & $0.03 \ |\ 4$      & $0.06 \ |\ 3$      & $0.08 \ |\ 3$     &  $0.23 \ |\ 5$  \\\hline
\end{tabular}
\end{table}
\vspace{-0.2cm}
\begin{table}[!ht]
 \centering
   \caption{Execution time (in seconds), and number of iterations. \\ Medium scale MDPs with $|\mathbb{A}| = 50$ actions. \label{tabB}}
\begin{tabular}{p{1.9cm}||c|c|c|c}
 & $|\mathbb{S}|=10^3$ & $|\mathbb{S}|=3.10^3$ & $|\mathbb{S}| = 5.10^3$ & $|\mathbb{S}|=8.10^3$   \\\hline \hline
RVI        & $8.17\ |\ 13839 $ & $24.33 \ |\ 17091 $ & $47.83 \ |\ 18026$ &  $71.72 \ |\ 17828 $  \\\hline \hline
 RPI + FP  & $277.23 \ |\ 6$   & $1004.23 \ |\ 6$    & $1818.46 \ |\ 6$   &  $2761.10 \ |\ 7 $ \\\hline
 RPI + GJ  & $14.58 \ |\ 6$    & $160.10 \ |\ 6$     & $585.87 \ |\ 6$    &  $2145.87 \ |\ 7$ \\\hline
 RPI + RB & $0.58 \ |\ 6$     & $2.04 \ |\ 6$       & $4.20 \ |\ 6$      &  $9.42 \ |\ 7$ \\\hline
\end{tabular}
\end{table}
\vspace{-0.2cm}
\begin{table}[!ht]
 \centering
   \caption{Execution time (in seconds), and number of iterations. \\
  Large scale MDPs with $|\mathbb{A}| = 100$ actions. \label{tabC}}
\begin{tabular}{p{1.9cm}||c|c|c|c}
  & $|\mathbb{S}|=10^4$ & $|\mathbb{S}|=5.10^4$ & $|\mathbb{S}| = 10^5$ & $|\mathbb{S}|=2.10^5$   \\\hline \hline
 RVI       & $125.27 \ |\ 17716 $ & $885.28 \ | \ 18876$  & $2862.14 \ | \ 20049$  & $6158.56 \ | \ 22079$  \\\hline \hline 
 RPI + FP  & $4700.23 \ |\ 8$     & $>10^4  \ |\ 8$       & $>10^4  \ |\ 9$       & $>>10^4  \ |\ 9$  \\\hline
 RPI + GJ  & $3517.49 \ |\ 8$     & $>10^4  \ |\ 8$       & $>10^4  \ |\ 9$       & $>>10^4  \ |\ 9$  \\\hline 
 RPI + RB  & $11.49 \ |\ 8$       & $183.70 \ |\ 8$        & $765.77 \ |\ 9$       & $3528.93 \ |\ 9$  \\\hline
\end{tabular}
\end{table}
\newpage
\subsection{Application : optimal policy analysis}
\begin{figure}[htbp]
    \centering
    \begin{minipage}{0.45\textwidth}
        \centering
        \includegraphics[width=\textwidth]{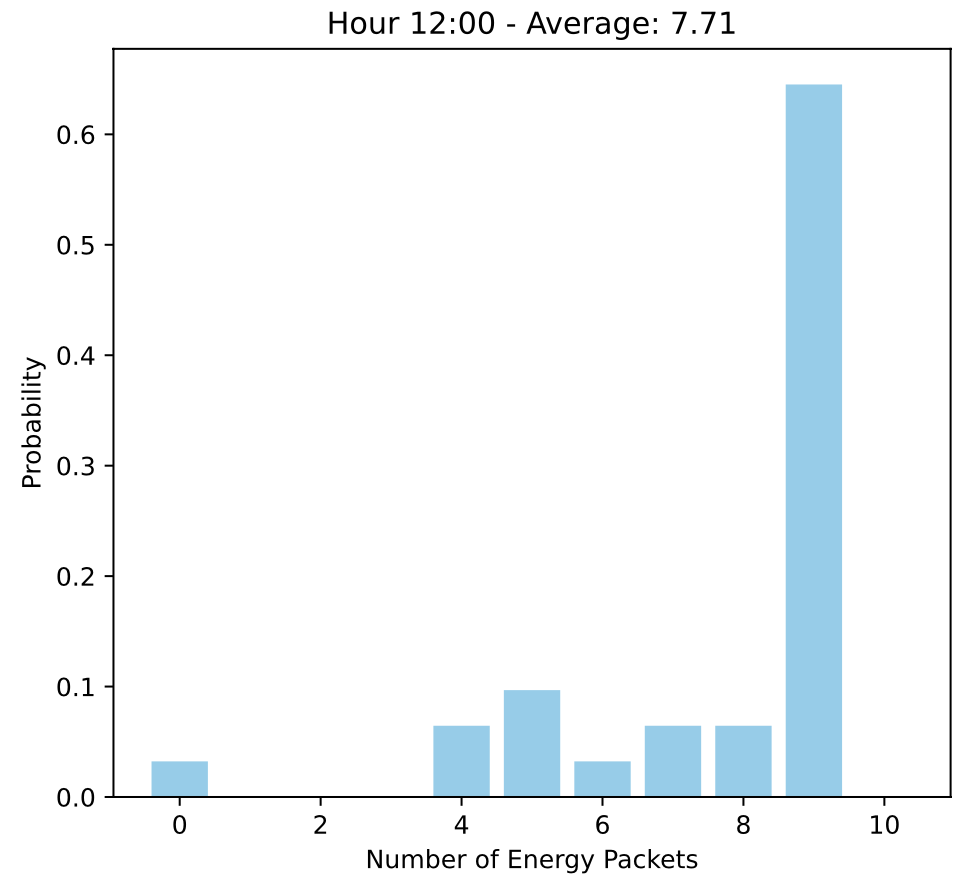}
    \end{minipage}
    \hfill
    \begin{minipage}{0.45\textwidth}
        \centering
        \includegraphics[width=\textwidth]{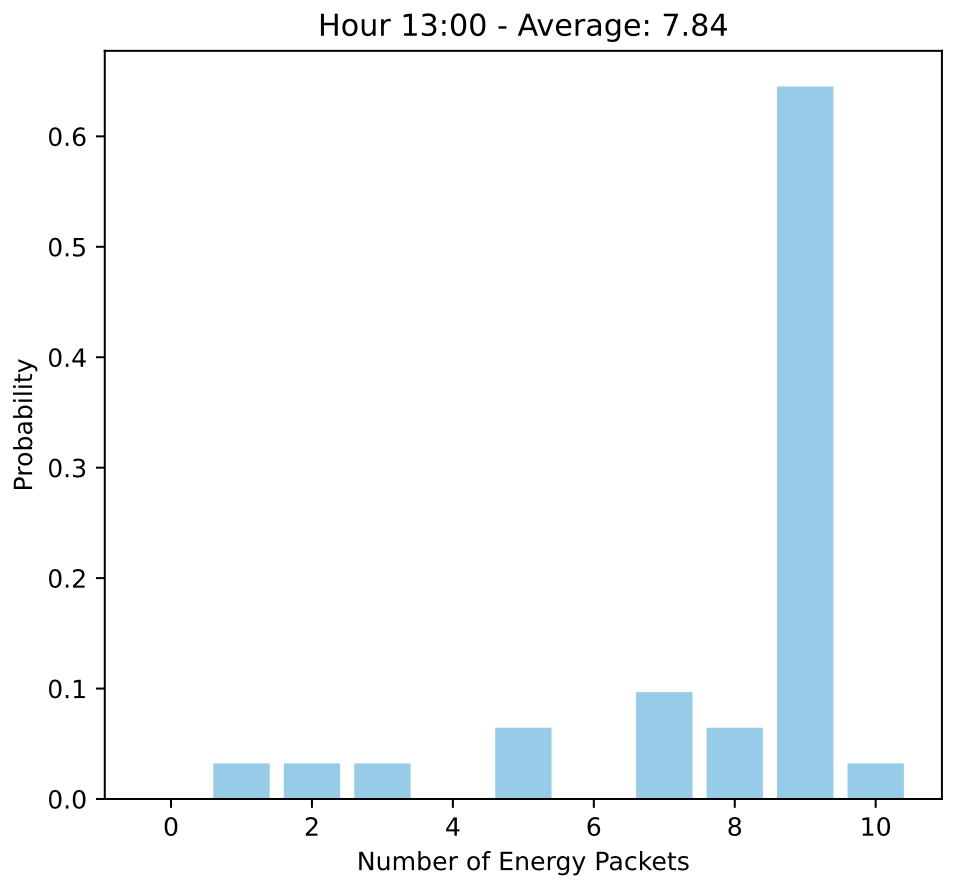}
    \end{minipage}
    \vspace{0.2cm}

    \begin{minipage}{0.47\textwidth}
        \centering
        \includegraphics[width=\textwidth]{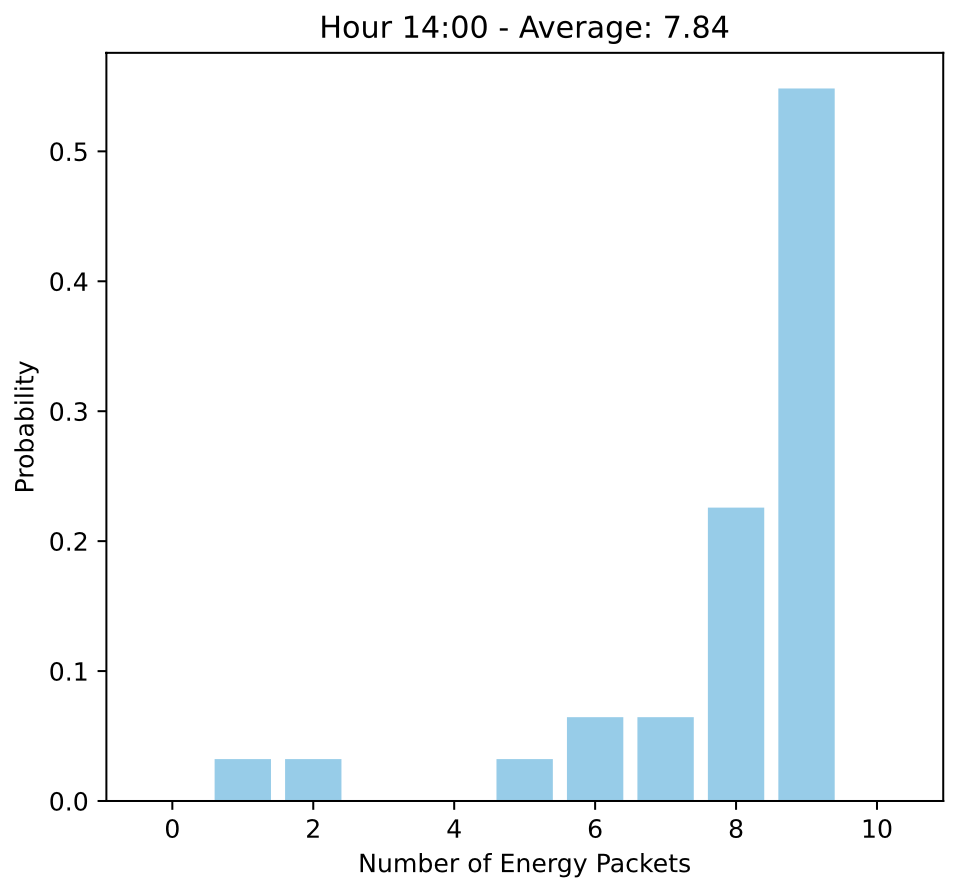}
    \end{minipage}
    \hfill
    \begin{minipage}{0.47\textwidth}
        \centering
        \includegraphics[width=\textwidth]{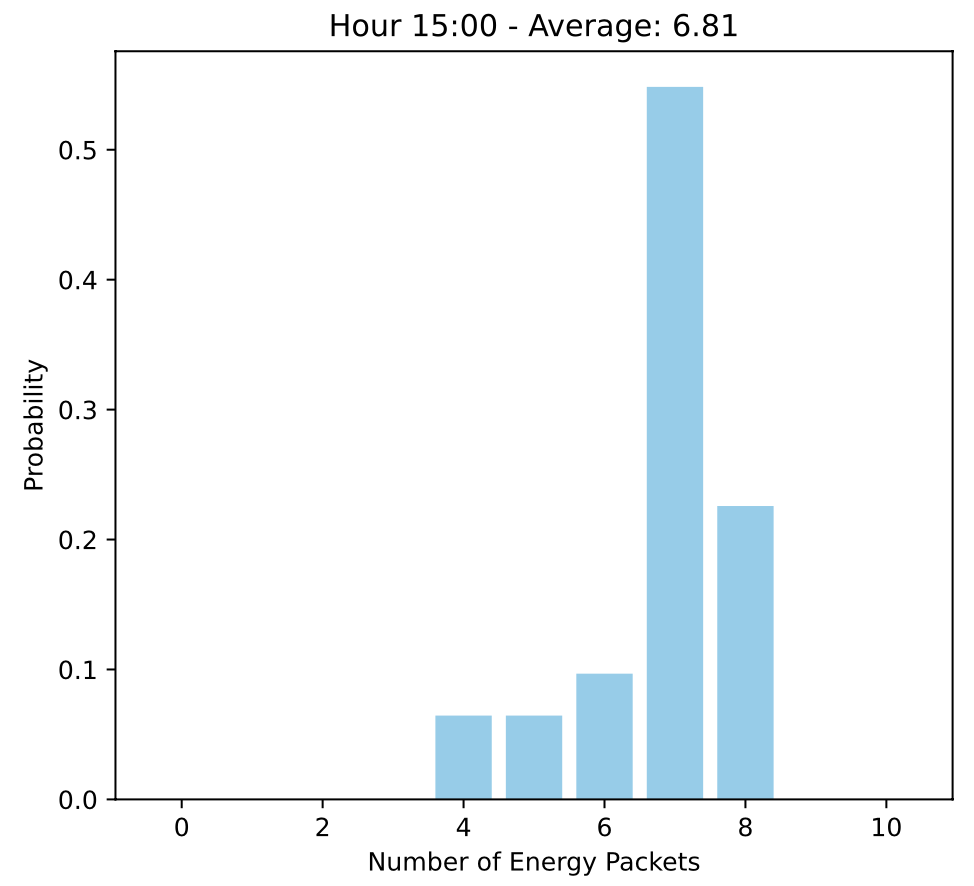}
    \end{minipage}
    \vspace{0.2cm}
    \begin{minipage}{0.45\textwidth}
        \centering
        \includegraphics[width=\textwidth]{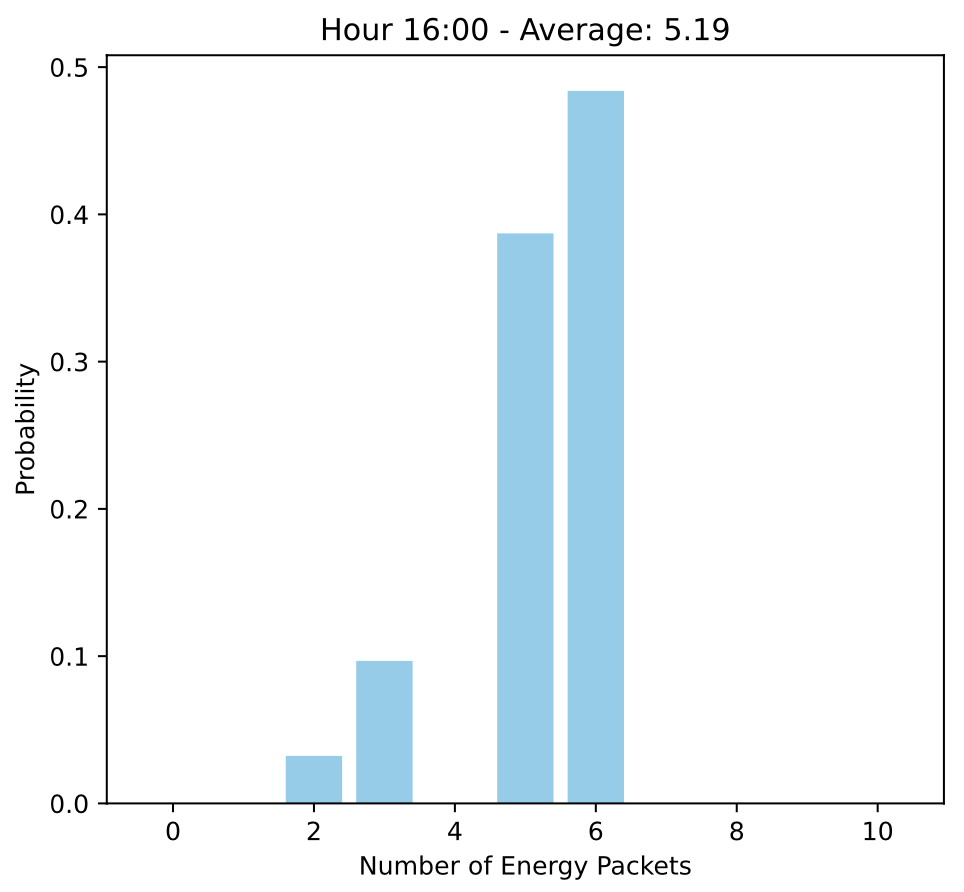}
    \end{minipage}
    \hfill
    \begin{minipage}{0.45\textwidth}
        \centering
        \includegraphics[width=\textwidth]{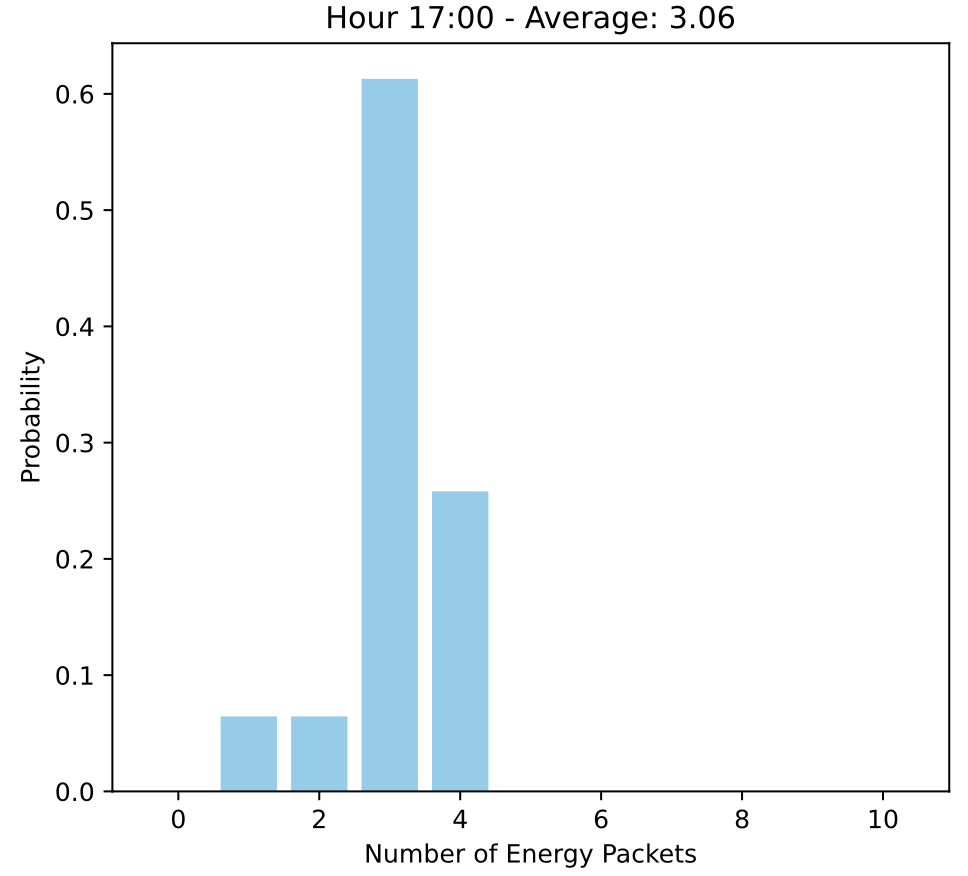}
    \end{minipage}
    \caption{Energy Packets distribution $\mathcal{A}_{12h}$, $\mathcal{A}_{13h}$ ... $\mathcal{A}_{17h}$. Barcelona, August.}
        \label{EP-Dists}
\end{figure}

To analyze the optimal behavior of the MDP agent, we base our approach on empirical data:
\begin{itemize}
\item Energy data preparation: we utilized open-access data from the PVWatts Calculator developed by the National Renewable Energy Laboratory (NREL \cite{NREL}), which provides hourly estimates of energy production in watt-hours (Wh) for a grid-connected solar panel at a specified location. This dataset incorporates various parameters, including climatic conditions sourced from the National Solar Radiation Database (NSRDB) and photovoltaic (PV) panel configurations, such as DC system size (kW), module type, system losses, azimuth angle, and tilt angle, among others. For data extraction, we retained the default PV system parameters provided by the calculator. Once the location is specified, hourly data can be downloaded. These data are based on a typical meteorological year file, representing a multi-year historical average for the set location and a fixed (open rack) photovoltaic system. From the downloaded data, we focused on the "AC System Output (W)" column, which indicates the power converted to alternating current (AC), suitable for network supply or storage in a battery. We therefore have data representing hourly energy production for each day and month over a typical year at a specific location. From these data we extract the following energy distributions.

    \item First, we conceptualized energy production in discrete energy packets, where each unit corresponding to "x" watt-hours. To calculate the number of energy units produced at a given hour, the total energy output at that time was divided by "x". Hence, this approach allows us to derive both the mean number of energy packets produced and the detailed distribution of energy packets generated for a specific hour and month across the dataset.
    An example of the average PV panel production in Barcelona in August is depicted in Fig. \ref{EP-arrival}. In this figure, the quantity of energy is represented in batches of $x=300$ Watts; for instance, at peak hour (14:00), the PV panel produces an average of $7.84$ energy packets, corresponding to $2352$ Wh, whereas the first energy quantity of the day (at 7:00) averages $0.9$ packets, or $270$ Wh. Notice that before 7:00 and after 18:00, the panel naturally no longer produces energy due to the absence of solar irradiations. 
    In this study, we analyze two types of curves related to energy packet production: (i) the first type represents the average number of energy packets produced per hour, providing a general overview of hourly production trends throughout the day; (ii) the second type represents detailed distributions of energy packets produced for each hour, capturing the variability in production at each specific hour of the day. For our MDP model, it is essential to consider these exact hourly distributions rather than average values, as they provide a finer understanding of production dynamics. Examples of these distributions are presented in Fig. \ref{EP-Dists}.

\item To model service demands, we use a distribution based on traffic intensity measured in Erlangs, specifically focusing on peak hours. We consider Fig. \ref{DP-arrival} distribution, inspired by empirical data from \cite{Cisco07} which reveals two distinct peaks in service demands occurring at 10:00 and 14:00.
    \end{itemize}
    
With the data prepared, we can now input it into our MDP model. For each hour, we have distinct distributions $\mathcal{A}_H$ representing energy packet production. For the month of August in Barcelona (with further comparisons across different locations and months provided in the next section), the time interval ranges from the first energy packet (EP) arrival at $t_0 = 7 \text{h}$ to the last arrival at $T = 18 \text{h}$, with energy packet productions described by distributions $\mathcal{A}_{7\text{h}}, \mathcal{A}_{8\text{h}}, \ldots, \mathcal{A}_{18\text{h}}$, some of which are illustrated in Fig. \ref{EP-Dists}. Similarly, service demands are represented by a Bernoulli probability drawn from $\mathcal{B}_H$, indicating the likelihood of a service request occurring. Other model parameters are set as follows: battery capacity $C = 65$, with a threshold level of $F = 25$ energy packets. The solar PV panel may experience failures due to manufacturing defects or premature wear. Based on real-world estimations, we assume a failure probability of $\approx 1.04\%$, resulting in transitions between PV-ON and PV-OFF states, modeled with probabilities $\alpha = 0.01$ (PV-ON to PV-OFF) and $\beta = 0.95$ (PV-OFF to PV-ON). We consider five possible actions associated with release probabilities $\mathcal{Z}_x^{M,(a)} = \{0.1, 0.3, 0.5, 0.7, 0.9\}$. For simplicity, these release probabilities are assumed to be consistent across both PV-ON and PV-OFF states. However, the filling process for energy packets differs, as we utilize an Interrupted Batch Process for arrivals. This configuration results in a model with $|\mathbb{S}| = 755$ states and $m = 4080$ arcs per action. The model is solvable efficiently, in less than one second, as shown in Tables \ref{tabA} and \ref{tabB}. Although this model is relatively small-scale, it is sufficient to analyze the agent’s behavior in optimizing rewards under various reward structures.

\begin{figure}[!p]
    \centering    \includegraphics[width=13cm,height=.43\textheight]{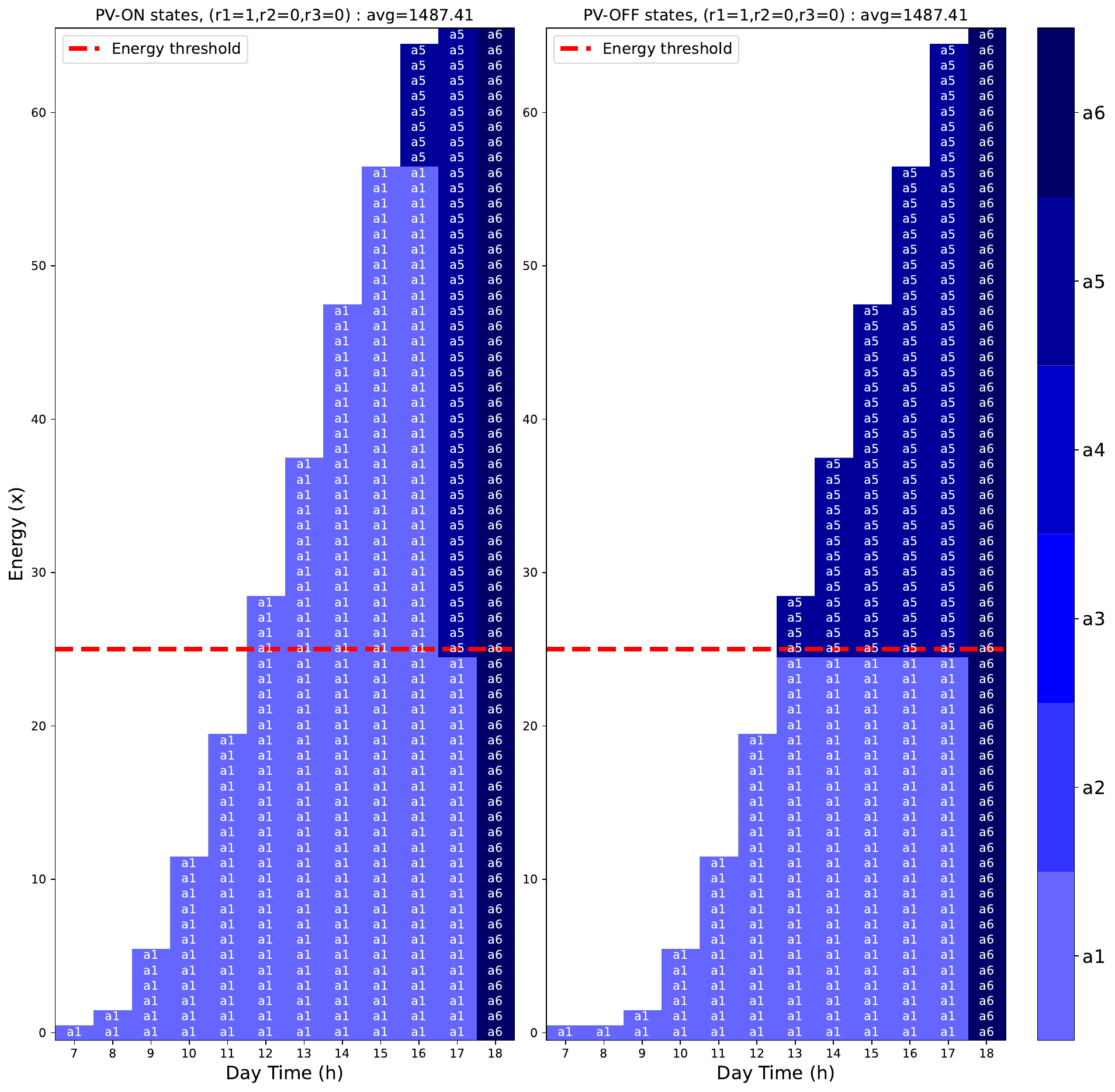}
    \vspace{-0.3cm}
    \caption{\label{Res1}Optimal policy of the system in separated heatmaps for PV-ON and PV-OFF states. Input rewards $r_1^+ = 1$, $r_2^- = 0$, $r_3^- = 0$.}
\end{figure}

\begin{figure}[!p]
    \centering
\includegraphics[width=13cm,height=.42\textheight]{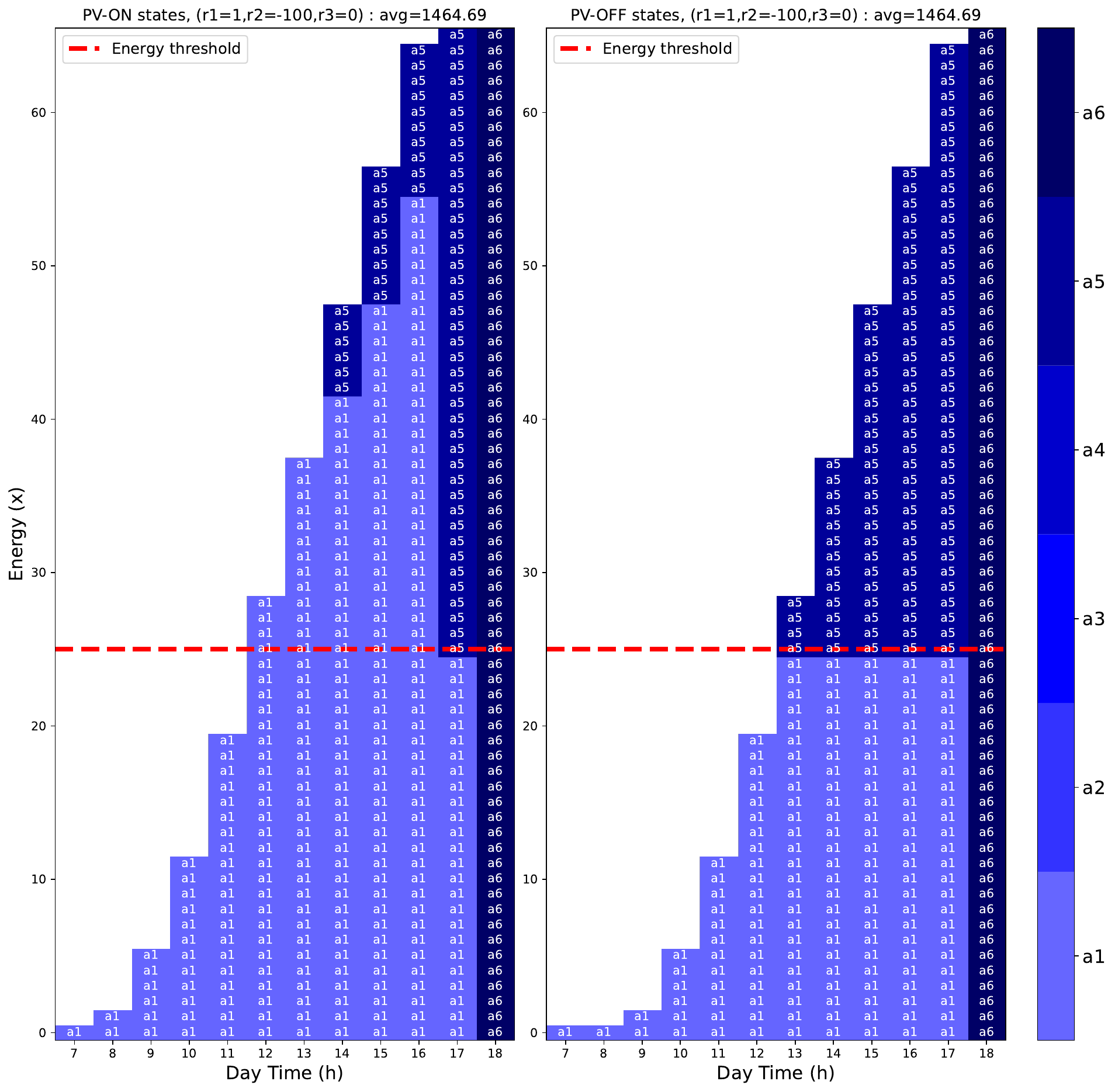}
    \vspace{-0.3cm}
    \caption{\label{Res2}Optimal policy of the system in separated heatmaps for PV-ON and PV-OFF states. Input rewards $r_1^+ = 1$, $r_2^- = -100$, $r_3^- = 0$.}
\end{figure}

\begin{figure}[!ht]
    \centering
\includegraphics[width=13cm,height=.42\textheight]{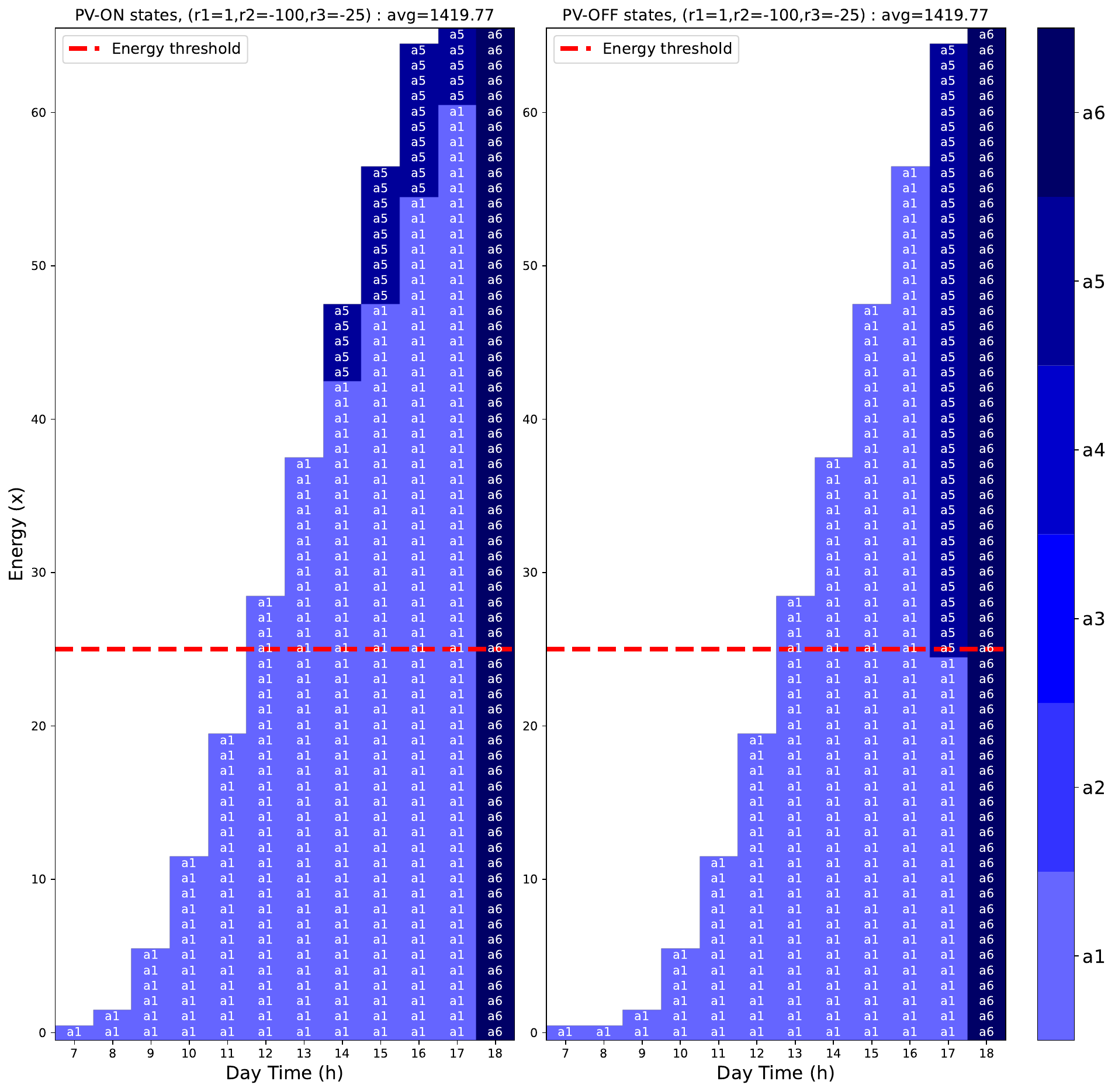}
    \caption{\label{Res3}Optimal policy of the system in separated heatmaps for PV-ON and PV-OFF states. Input rewards : $r_1^+ = 1$, $r_2^- = -100$, $r_3^- = -25$.}
\end{figure}

In following experiments, we first (Figure \ref{Res1}) consider only reward $r_1^+ = 1$ (a positive reward for battery release, i.e., battery sold). The second experiment, (Figure \ref{Res2}), incorporates both $r_1^+ = 1$ and $r_2^- = -100$ (a negative reward for energy packet loss), and the third experiment (Figure \ref{Res3}) additionally includes $r_3^- = -50$ penalizing the agent when no energy packets are present in the battery to serve an actual task arriving to the network. Each figure presents two heatmaps illustrating the optimal policy for each experiment. The left heatmap corresponds to the system during PV-ON states, and the right heatmap refers to PV-OFF states. From these experiments, we observe the following.

\begin{itemize}
\item The agent exhibits distinct behaviors depending on whether the system is operating during PV-ON or PV-OFF states. Specifically, during PV-ON, the agent tends to wait until the battery is more fully charged before initiating a release. Conversely, during PV-OFF, in the absence of energy packet arrivals, the decision shifts towards releasing the battery as soon as it reaches threshold $F=25$, indicated by the dotted red line. This behavior, consistently observed across all three experiments, is predominantly influenced by the reward $r_1^+$, which incentivizes the agent to release a well-charged battery for sale. 
For instance, in non-failure states of Experiment 1, the optimal decision is to release the battery more frequently in the later hours of the day, specifically when $h = 16$ and $x \geq 57$, or $h=17$. Additionally, when $h = T = 18$, the end-of-day deadline is reached, resulting in an automatic release of the battery (denoted by dark blue "a6" in the heatmaps). In Experiment 2, however, the release decision occurs earlier in the day at $h \geq 14$ and $x \geq 43$, or $h=15$ with $x \geq 48$, or $h=16$ with $x \geq 55$, and again at $h=17$. This shift is attributed to the penalty $r_2^-$, which emphasizes the importance of preventing packet losses, leading the agent to release the battery earlier in the day. This strategy ensures that packets, which would otherwise be lost, are accommodated by the replacement battery after the initial one is released, accounting for the prevalence of dark blue in PV-ON heatmap of Experiment 2. 
It is also noteworthy that no differences are observed in PV-OFF states between Experiments 1 and 2, as the penalty $r_2^-$ is triggered solely by packet arrivals that could lead to losses, while during PV-OFF states, there are no arrivals due to PV failure. In the final experiment, we introduce the penalty $r_3^-$, which applies when the battery is empty, potentially discouraging service delays. This penalty influences the agent’s decisions by discouraging premature battery releases, aiming to avoid penalties associated with an empty battery and ensuring timely data packets service. Consequently, the optimal release threshold is adjusted to $h \geq 14$ and $x \geq 45$, or $h=15$ with $x \geq 48$, or $h=16$ with $x \geq 55$, or $h=17$ with $x \geq 61$. In PV-OFF states, the agent maintains the release only at $h=17$. In sum, the agent has reduced the number of releasing states (i.e., dark blue areas) in response to this new penalty.
Those adjustments reflects the agent's strategy to balance the need to release batteries against the risk of leaving the battery empty which would encourages data packets delays, and the risk of overfilling, which could result in penalties due to energy packet losses. This approach effectively reduces the likelihood of incurring penalties while maintaining adequate gains.

\item One can also observe the value of $\rho^{(\pi^*)}$, the average reward, displayed above the figures. This value decreases progressively across the experiments, from $1487.41$ in Experiment 1, to $1464.69$ in Experiment 2, and further to $1419.77$ in Experiment 3. This decrease is provoked by the various penalties introduced in each successive experiment. $\rho^{(\pi^*)}$ is also related to the specific magnitudes selected for each immediate reward; typically, we have set $r_1=1$, $r_2=-100$, and $r_3=-25$. This selection is designed to more effectively illustrate the trade-off act between these three rewards in the determination of the optimal policy. 
\end{itemize}


\newpage
\subsection{Application : cities comparison}
Next, we aim to analyze which location and month of the year would be most suitable for an Off-Grid telecommunication operator to deploy its infrastructure. We assume the same non-stationary distribution for service demands (Fig. \ref{DP-arrival}) but allow for varying energy production distributions depending on location and month. The immediate rewards are set to $r_1^+ = 1$, $r_2^- = -100$, and $r_3^- = -200$, introducing a slightly stronger emphasis on data packet delays compared to previous experiments. 

\begin{figure}[p]
\vspace{-2cm}
    \centering
    \begin{subfigure}{0.9\textwidth}
        \centering
    \includegraphics[width=0.9\linewidth]{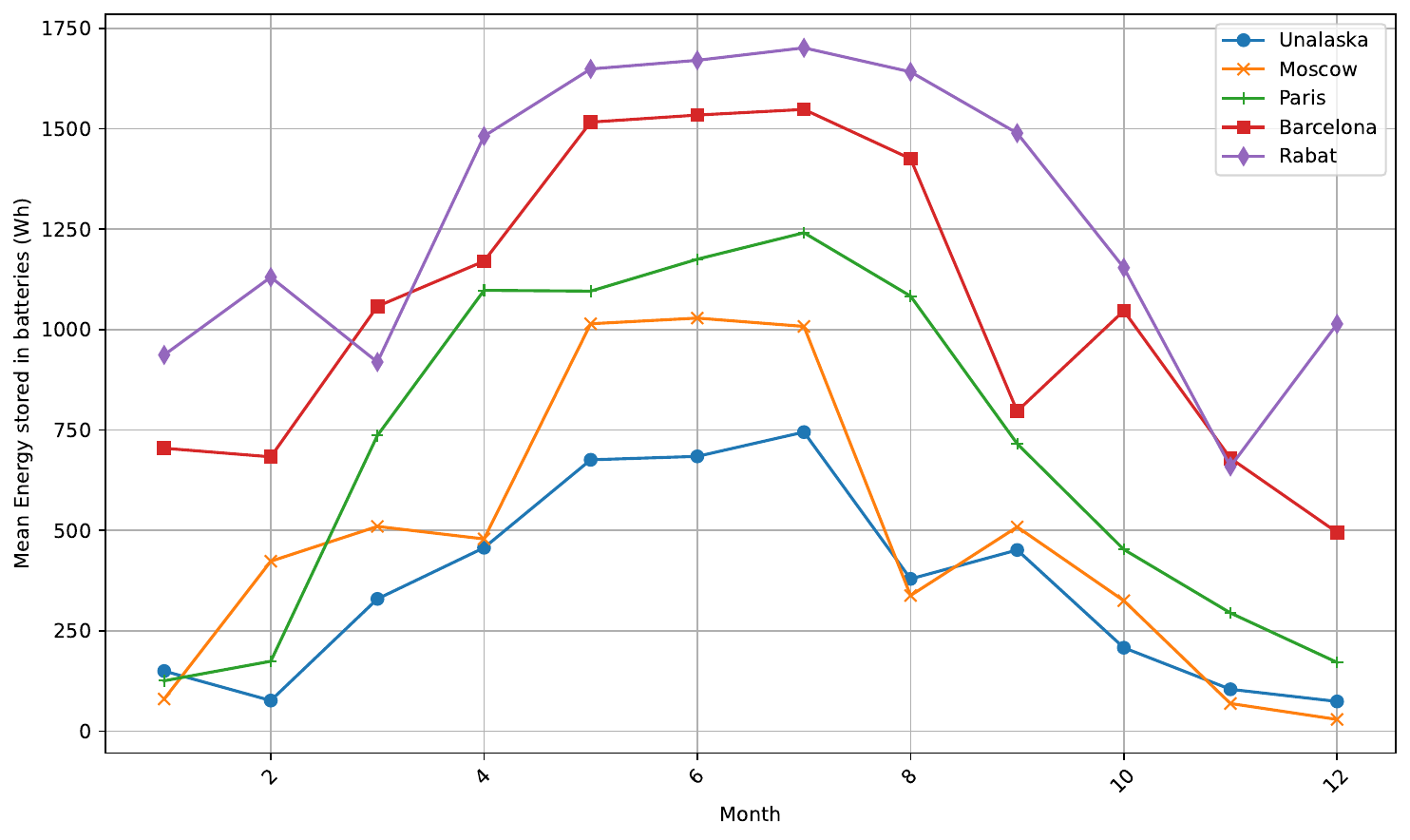}
    \vspace{-0.2cm}
        \caption{Average quantity of energy stored in batteries at release, $\mathbb{E}[Release]$}
        \label{fig:Erelease}
    \end{subfigure}
    \hfill
    \begin{subfigure}{0.9\textwidth}
        \centering
        \includegraphics[width=0.9\linewidth]{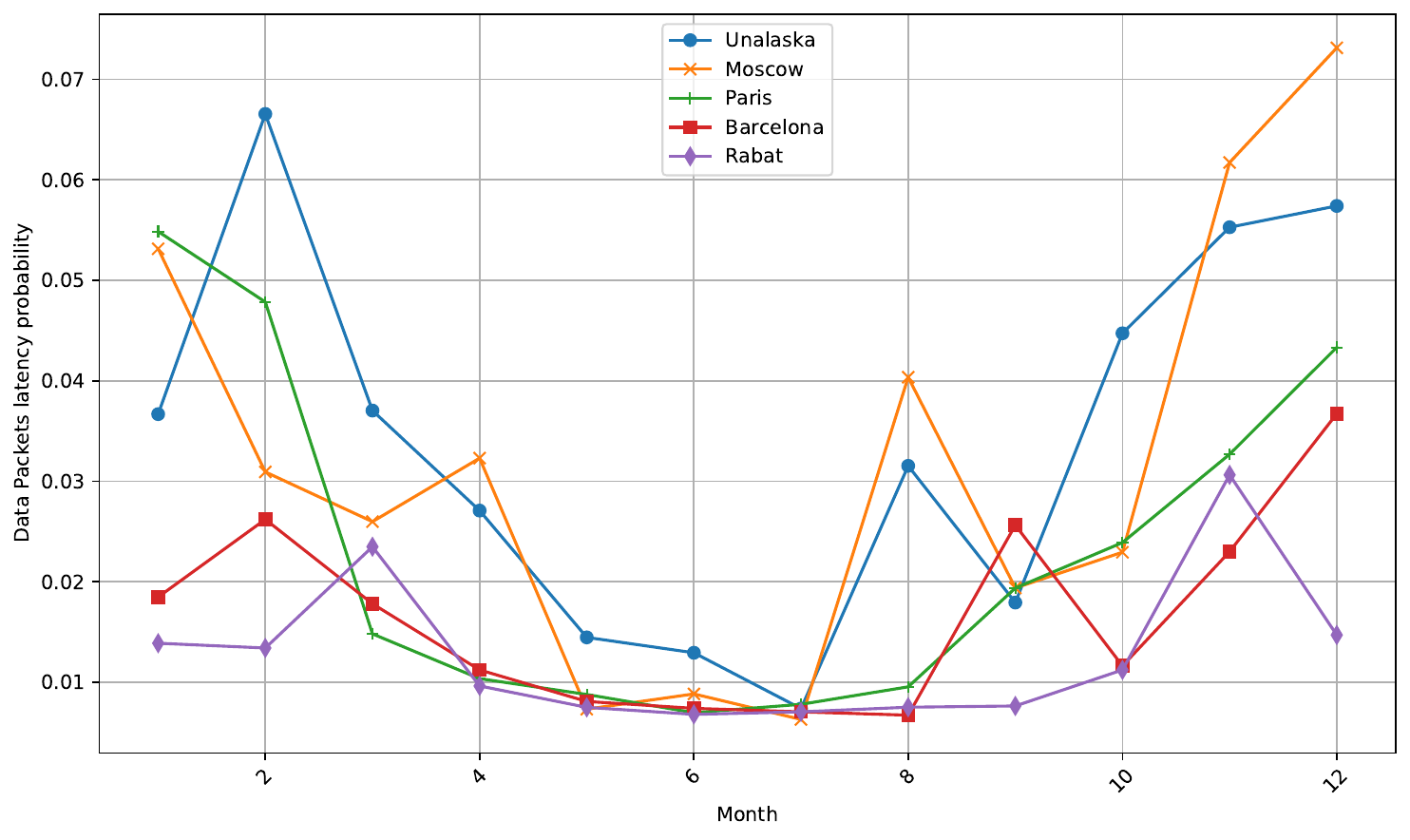}
            \vspace{-0.2cm}
        \caption{Delay probability, $\mathbb{E}[Delay]$}
        \label{fig:Edelay}
    \end{subfigure}
    \begin{subfigure}{0.9\textwidth}
        \centering
        \includegraphics[width=0.9\linewidth]{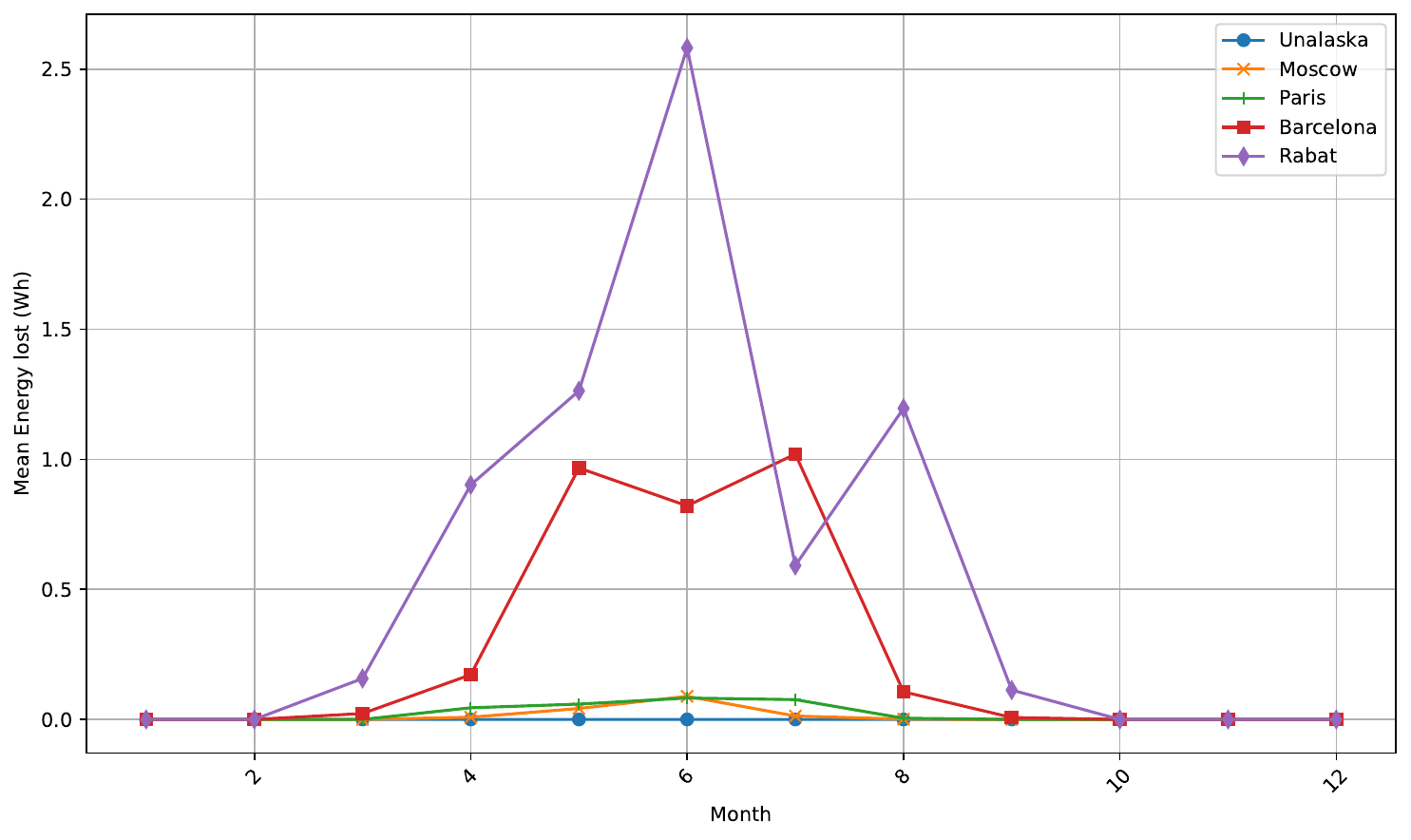}
        \vspace{-0.2cm}
        \caption{Average quantity of energy lost, $\mathbb{E}[loss]$}
        \label{fig:Eloss}
    \end{subfigure}
    \vspace{-0.2cm}
    \caption{Optimal measures across different locations throughout a typical year}
    \label{fig:measures}
\end{figure}

\begin{figure}[!ht]
    \centering
    \includegraphics[width=0.8\linewidth]{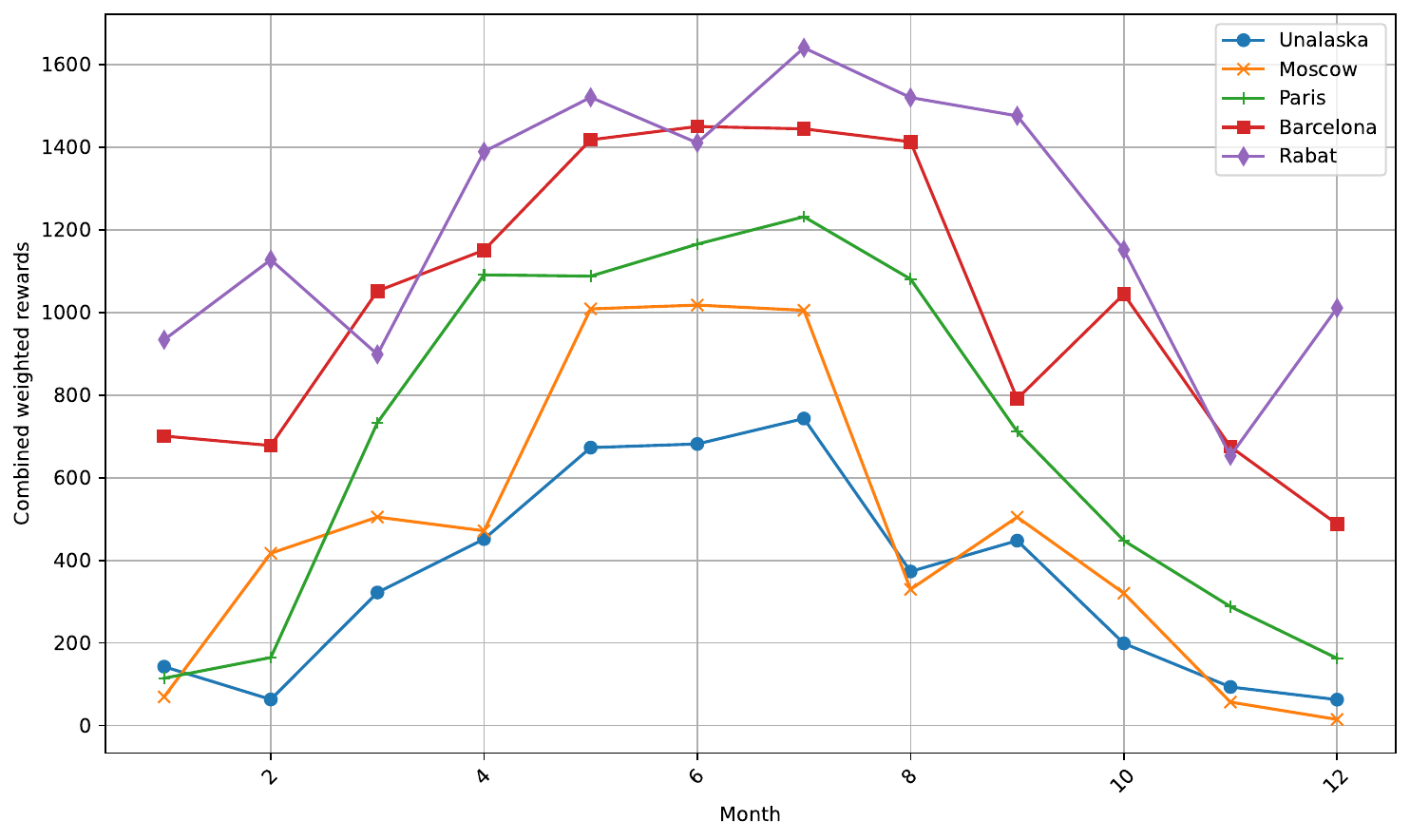}
    \caption{Combined optimal measures across different locations throughout a typical year}
    \label{fig:combined}
\end{figure}

In Fig. \ref{fig:measures}, we compare measures for various cities: Rabat, Barcelona, Moscow, Paris, and Unalaska. These metrics, as defined in equations \eqref{Erelease}, \eqref{ELost}, and \eqref{EDelay}, represent the average amount of energy stored in batteries (in watt-hours), the probability of delay, and the average amount of energy loss (in watt-hours), respectively. 
We first observe the impact of geographical location and climatic conditions on battery storage. For instance, Unalaska experiences lower solar irradiance compared to Rabat or Barcelona, as clearly illustrated in Batteries filling in Fig. \ref{fig:Erelease}. Additionally, battery storage levels peak during the summer months, indicating higher energy availability and lower system delays during this period, as shown in Fig. \ref{fig:Edelay}. However, energy losses are also most pronounced in this period of the year, particularly in Rabat and Barcelona, as depicted in Fig. \ref{fig:Eloss}.
A telecommunications operator might also consider a combination of these three metrics, illustrated in Fig. \ref{fig:combined} where the aim is to maximize the combined rewards. Here, Rabat and Barcelona demonstrate the best results, primarily driven by the higher average energy stored in batteries that comes from solar irradiance. However, this behavior is not stationary, for instance Rabat shows mostly the best results along the year except in March and June which can be related to energy losses during these months, with a peak of energy loss in June.


\newpage
The proposed MDP analysis, datasets, source code and results are available in GitHub \cite{YAEM24}.

\section{Conclusion}

In this study, we show that utilizing Markov chain structures offers a practical method for creating efficient decision-making solutions, especially in telecommunications. As simpler telecom nodes controlled by intelligent agents face significant time and energy limitations, using the natural structure of these problems allows for the development of both energy-efficient and computationally effective algorithms. In this work, we based our approach on non-stationnary arrivals with failure states to match real conditions while using empirical data. We have proposed a complete methodology from data preparation, efficient model-solving to decision making.
Future work will extend this methodology by investigating Robertazzi's type A structure \cite{Rob90} associated with Hessenberg matrices (see \cite{Stew94}), tensor decomposition of chains \cite{PFLe88}, and chain aggregation based on graph properties \cite{Stew94}. To the best of our knowledge, such strategies have not yet been applied to the analysis of decision problems. In this context, we aim to compare Off-Grid systems with hybrid models \cite{GWBO24} that integrate grid connectivity, analyzing the decision-making process between relying on local supply from PV panels, subject to meteorological conditions, and grid supply, leading to additional costs.

\newpage
\bibliographystyle{elsarticle-num} 
\bibliography{ComCom2025_Version_MDP_Battery}

\end{document}